\documentclass[final]{siamltex1213}

\usepackage{xspace,exscale}
\usepackage{fancybox}
\usepackage{graphicx}
\usepackage{color}
\usepackage{amsfonts}
\usepackage{amssymb}
\usepackage{url}

\usepackage[noadjust]{cite}

\usepackage{amsmath}
	\makeatletter
	\let\over=\@@over \let\overwithdelims=\@@overwithdelims
	\let\atop=\@@atop \let\atopwithdelims=\@@atopwithdelims
  	\let\above=\@@above \let\abovewithdelims=\@@abovewithdelims
  	\makeatother
\usepackage{rotating}

\usepackage{ifpdf}

\usepackage{subfigure}
\usepackage{psfrag}

\newcommand{\matc}{\ensuremath{\mathcal{C}}}
\newcommand{\matx}{\ensuremath{\mathcal{X}}}

\newcommand{\matn}{\ensuremath{\mathcal{N}}}

\newcommand{\mcomplex}{\ensuremath{\mathbb{C}}}

\newcommand{\FF}{\ensuremath{\mathbb{F}}}
\newcommand{\Sph}{\ensuremath{\mathbb{S}}}

\ifx\eqref\undefined
	\newcommand{\eqref}[1]{~(\ref{#1})}
\fi
\ifx\mod\undefined
	\def\mod{\mathop{\rm mod}}
\fi

\def\dperp{\perp\!\!\!\perp}
\newcommand{\Bern}{\text{Bern}}

\def\dim{\mathop{\rm dim}}

\def\exp{\mathop{\rm exp}}

\def\sign{\mathop{\rm sgn}}

\def\EE{\mathbb{E}\,}

\def\PP{\mathbb{P}}

\def\eqdef{\stackrel{\triangle}{=}}

\def\unifto{\mathop{{\mskip 3mu plus 2mu minus 1mu%
	\setbox0=\hbox{$\mathchar"3221$}%
	\raise.6ex\copy0\kern-\wd0%
	\lower0.5ex\hbox{$\mathchar"3221$}}\mskip 3mu plus 2mu minus 1mu}}

\ifx\lesssim\undefined
\def\simleq{{{\mskip 3mu plus 2mu minus 1mu%
	\setbox0=\hbox{$\mathchar"013C$}%
	\raise.2ex\copy0\kern-\wd0%
	\lower0.9ex\hbox{$\mathchar"0218$}}\mskip 3mu plus 2mu minus 1mu}}
\else
\def\simleq{\lesssim}
\fi

\ifx\gtrsim\undefined
\def\simgeq{{{\mskip 3mu plus 2mu minus 1mu%
	\setbox0=\hbox{$\mathchar"013E$}%
	\raise.2ex\copy0\kern-\wd0%
	\lower0.9ex\hbox{$\mathchar"0218$}}\mskip 3mu plus 2mu minus 1mu}}
\else
\def\simgeq{\gtrsim}
\fi

\newif\ifmapx
{\catcode`/=0 \catcode`\\=12/gdef/mkillslash\#1{#1}}
\edef\jobnametmp{\expandafter\string\csname hc_hamming_apx\endcsname}
\edef\jobnameapx{\expandafter\mkillslash\jobnametmp}
\edef\jobnameexpand{\jobname}
\ifx\jobnameexpand\jobnameapx
\mapxtrue
\else
\mapxfalse
\fi

\long\def\apxonly#1{\ifmapx{\color{blue}#1}\fi}

\begin{document}

\title{Hypercontractivity of spherical averages in Hamming space}

\author{Yury Polyanskiy\footnotemark[2]}

\maketitle
\renewcommand{\thefootnote}{\fnsymbol{footnote}}
\footnotetext[2]{YP is with the Department of Electrical Engineering 
and Computer Science, MIT, Cambridge, MA 02139 USA.
	\mbox{e-mail:~{\ttfamily yp@mit.edu}}.\\
The research was supported by the NSF grant CCF-12-53205 and NSF Center for Science of Information (CSoI) 
under grant agreement CCF-09-39370.}
\renewcommand{\thefootnote}{\arabic{footnote}}

\begin{abstract}
Consider the linear space of functions on the binary hypercube and the linear operator $S_\delta$ acting
by averaging a function over a Hamming sphere of radius $\delta n$ around every point. It is shown that this
operator has a dimension-independent bound on the norm $L_p \to L_2$ with $p = 1+(1-2\delta)^2$. This
result evidently parallels a classical estimate of Bonami and Gross for $L_p \to L_q$ norms for the
operator of convolution with a Bernoulli noise. The
estimate for $S_\delta$ is harder to obtain since the latter is neither a part of a semigroup, nor a tensor power. 
The result is shown by a detailed study of the eigenvalues of $S_\delta$ and $L_p\to L_2$ norms of the
Fourier multiplier operators $\Pi_a$ with symbol equal to a characteristic function of the Hamming
sphere of radius $a$ (in the notation common in boolean analysis $\Pi_a f=f^{=a}$, where $f^{=a}$ is a degree-$a$
component of function $f$).
A sample application of the result is given: Any set $A\subset \FF_2^n$ with the property that
$A+A$ contains a large portion of some Hamming sphere (counted with multiplicity)
must have cardinality a constant multiple of $2^n$. 
\end{abstract}



\section{Main result and discussion}

Consider a linear space $\mathcal{L}$ of functions on $n$-dimensional
Hamming cube $f:\FF_2^n\to\mcomplex$. We endow $\mathcal{L}$ with the following norms and an inner product:
\begin{align} \|f\|_p &\eqdef \EE^{1\over p}[ |f(X)|^p]\,, \qquad
1\le p \le \infty\,,\\
	 (f, g) &\eqdef \EE[f(X) \bar g(X)]\,,
\end{align}
where $X$ is uniform on $\FF_2^n$.
For any linear operator $T:\mathcal{L}\to\mathcal{L}$ we define
	$$ \|T\|_{p\to q} \eqdef \sup_{f\in \mathcal{L}} {\|Tf\|_q\over \|f\|_p}\,.$$
\apxonly{(These are not Banach-algebra norms unless $p=q$.)}
Let $Z=(Z_1,\ldots,Z_n)$ be a random element of $\FF_2^n$ with components independent and identically distributed (i.i.d.)
according to $\Bern(\delta)$ distribution: $\PP[Z_i=1]=1-\PP[Z_i=0]=\delta$. For the following operator
\begin{equation}\label{eq:ndelta}
	 N_\delta f(x) \eqdef \EE[f(x+Z)], \qquad x \in
\FF_2^n\,, 0\le \delta\le 1
\end{equation}
the so-called ``hypercontractive'' inequality was established by Bonami~\cite{AB70}, Gross~\cite{LG75} and
others (see~\cite[Chapter 9, notes]{odonnell-book} for the history):
\begin{equation}\label{eq:bonami}
	 \|N_\delta f\|_{q} \le \|f\|_p \,, \qquad \forall q \ge p \ge 1, p -1 \ge (q-1)(1-2\delta)^2, p,q\ge1.
\end{equation}

There are a number of applications of hypercontractive inequalities. For example, we mention an early result in
information theory~\cite{AC76}, which has recently become known as the ``small-set expansion''. A number of applications in theoretical computer
science are presented in~\cite[Chapter 9-10]{odonnell-book}. One of the pillars of the analysis of boolean functions, the KKL
lemma~\cite{KKL88}, is an ingenious application of~\eqref{eq:bonami}. Hypercontractivity is also an indispensable tool
in probability for analyzing mixing of Markov chains~\cite{DSC96} and isoperimetry~\cite[Theorem 3.4]{MORSS-nicd}.

In this paper we analyze the $L_p \to L_2$ norm for an operator $S_{\delta}$ of averaging over a Hamming
sphere $\Sph_{\delta n}$. Specifically, for $x=(x_1,\ldots,x_n)\in \FF_2^n$ denote the Hamming weight of
$x$ and the Hamming sphere centered at zero as
\begin{align} |x| &\eqdef |\{j: x_j=1\}|\\
	 \Sph_j &\eqdef \{x: |x|=j\} \,.
\end{align}
The operator $S_{\delta}$ is defined as follows:
	$$ S_\delta f(x) \eqdef {n \choose j }^{-1}
				\sum_{y\in\FF_2^n, |y| = j } f(x+y)\,,$$
where $j=\lceil \delta n\rceil$ if $\delta <1/2$ and $j=\lfloor \delta n\rfloor$ if $\delta \ge 1/2$.
In other words, we may write
	$$ S_\delta f \eqdef {f*1_{\Sph_j} \over |\Sph_j|}\,,$$
where $*$ denotes the convolution
	$$ f*g(x) \eqdef \sum_{y\in\FF_2^n} f(x-y) g(y)\,. $$
This definition ensures $S_\delta f (x) = S_{1-\delta} f (\bar x)$ for $\delta \neq {1\over 2}$, where $\bar
x=(1-x_1,\ldots,1-x_n)$.

Our main result is that $S_\delta$ satisfies an inequality entirely similar to $N_\delta$, namely:
\begin{equation}\label{eq:new_hc_eq}
	 \|S_\delta f\|_{2} \le C_\delta \|f\|_p \,, \qquad \forall  p \ge 1+(1-2\delta)^2, \delta \neq {1\over 2}\,,
\end{equation}
where the crucial part is that $C_\delta>1$ does not depend on dimension $n$. Note also that the constant cannot be tightened to 1. Indeed, taking $f = 1_{even}$ to be the characteristic function
of the set of all even-weight vectors we get  
	$$ \|S_\delta\|_{p\to 2} \ge 2^{{1\over 2} - {1\over p}}\,,\qquad 1\le p \le 2\,, 0<
\delta <1\,, $$
regardless of dimension $n$. More precisely, we show the following.

\begin{theorem}\label{th:main} Consider the set $F \subset [0,1]\times[1,2]$
	$$ F = \{(\delta, p): p\ge 1 + (1-2\delta)^2, 0\le \delta \le 1, 1 < p \le 2\}\,.$$
	For every compact subset $K$ of $F$ there exists a constant $C=C(K)$ such that for all
$(\delta,p)\in K$, $n\ge 1$ and $f: \FF_2^n \to \mcomplex$ we have
\begin{equation}\label{eq:main}
	 \| S_\delta f\|_2 \le C \|f\|_p\,.
\end{equation}
Conversely, for any $(\delta,p) \not\in F$ there is $E>0$ such that
\begin{equation}\label{eq:main_cex}
	 \sup_f {\| S_\delta f\|_2  \over \|f\|_p} \ge e^{n E +  o(n)}\,,\qquad n\to\infty
\end{equation}
with the exception of $\delta=1/2, p=1$ for which we have
\begin{equation}\label{eq:main_young}
	 \sup_f {\| S_{1/2} f\|_2  \over \|f\|_1}  = 2^{n/2} {n \choose \lfloor n/2\rfloor}^{-{1\over2}} \sim
\left(\pi n\over 2\right)^{1\over 4}\,.
\end{equation}
\end{theorem}
\textbf{Remark:} The constants that can be extracted from our proof method (after numerical evaluations) are as follows: for $\delta\le 0.16$ we
have $C = \sqrt{2}$, while for larger $\delta$ we can take $C$ to be arbitrarily close to $\sqrt{2}$
for sufficiently large $n$.

The full proof is given in Section~\ref{sec:proof}, while here we give a high-level sketch.  We note first that the standard methods for showing hypercontractivity do not apply since they require the operator to be a tensor
power or be part of a semigroup. The semigroup could be continuous-time, as in~\cite{DSC96}, or discrete-time as
in~\cite{LM97}, but $S_\delta$ is a member of neither. \apxonly{Orig: 
that our operators $S_\delta$ do not form a semigroup. This may potentially be worked around by
applying the discrete-time version of log-Sobolev inequalities developed by Miclo~\cite{LM97}.
However, the natural comparison to $N_\delta$ via Miclo's method is unfortunately not useful: the
primary reason is that the log-Sobolev constant of $N_\delta$ is of order ${1\over n}$ which implies
tight hypercontractive estimates when $\delta \sim {1\over n}$ and is very loose otherwise. -- I have no idea what I
mean by log-Sobolev constant is of order $1\over n$.\\\par} Instead, our proof proceeds by noticing that $S_\delta$
and $N_\delta$ commute and are self-adjoint, hence have common orthogonal eigenspaces (given by the Fourier transform,
also known as degree-$d$ components).
Consequently, decomposing a function $f=\sum_j f_j$ into sum of its projections on eigenspaces we have
from~\eqref{eq:bonami}:
\begin{equation}\label{eq:rfe0}
	\|N_\delta f\|_2^2 = \sum_j \lambda_j(N_\delta)^2 \|f_j\|_2^2 \le \|f\|_p^2\,.
\end{equation}
Writing a similar expansion for $S_\delta$ we have
\begin{equation}\label{eq:rfe1}
	\|S_\delta f\|_2^2 = \sum_j \lambda_j(S_\delta)^2 \|f_j\|_2^2\,.
\end{equation}
If we had that $\lambda_j(S_\delta) \le \lambda_j(N_\delta)$, then we could just upper bound~\eqref{eq:rfe1}
with~\eqref{eq:rfe0} and conclude the proof. It turns out that such estimate does hold but only for a range of $j$, and
thus the bulk of the proof consists of showing that contribution to~\eqref{eq:rfe1} of the eigenspaces outside of this
range is small. This part crucially depends on a curious relation between
norms of certain Fourier-multiplier operators on $\FF_2^n$ and eigenvalues of $S_\delta$. The corresponding
estimates that bound energies in the degree-$a$ components of functions on the hypercube are, perhaps, of independent
interest.

\subsection{Discussion}
Why would one conjecture that $S_\delta$
is hypercontractive? Note that~\cite[Theorem 3.7]{DSC96} shows that a discrete time Markov
chain on state space $\matx$ and whose kernel satisfies hypercontractive inequality, mixes in time of 
order $O(\log\log |\matx|)$. For $S_\delta$, this Markov chain is a non-standard random walk on a
hypercube $\FF_2^n$ which jumps by a distance exactly $\delta n$ at each step. A simple coupling
argument shows that indeed such a random walk must mix in time $O(\log n)$, therefore giving some 
probabilistic intuition as to why Theorem~\ref{th:main} might hold. 

We note that our main goal was to show an $O(1)$ estimate for $\|S_\delta\|_{p\to q}$. Indeed, a $O(\sqrt{n})$ estimate
is much easier:
\begin{theorem} For any $\delta$ and $p \ge 1 + (q-1)(1-2\delta)^2$ we have
		$$ \|S_\delta\|_{p\to q} = O(\sqrt{n})\,.$$
\end{theorem}
\begin{proof} Assuming without loss of generality that $f\ge 0$ it is easy to see from Stirling's formula that
	$$ {1\over {n\choose \delta n}} \sum_{|y|=\delta n} f(x+y) \le O(\sqrt{n}) \sum_{|y|=\delta
n} f(x+y) \delta^{|y|} (1-\delta)^{n-|y|}\,.$$
	Then extending summation to all of $y$ we get 
	$$ S_\delta f (x) \le O(\sqrt{n}) N_\delta f (x) \qquad \forall x\in\FF_2^n\,.$$
	The result then follows from~\eqref{eq:bonami}.
\end{proof}

The importance of having an $O(1)$ estimate for the $p\to q$ norm is due to the following general result of 
Semenov and Shneiberg~\cite{SS88}, which generalized earlier results of Fefferman and Segal~\cite{FS72,IS70}. \apxonly{
Note that one of the most fascinating properties of~\eqref{eq:bonami} is that it shows the following
``stickiness at 1'' of $\|\cdot\|_{p\to q}$ norms: as operator $N_\delta$ starts to depart from the
$N_{1\over 2}$ the norm $\|N_\delta\|_{p\to q}$ remains stuck at 1 for some range of values
$\delta\in[\delta_0, {1\over2}]$ before starting to
grow as $\delta <\delta_0$. This distinguishes the measure of dependence $\|\cdot\|_{p\to q}$ from other measures (such as
mutual information, or correlation coefficients). Interestingly, a similar effect was observed for
Fourier multiplier operators and norms $\|\cdot\|_{p\to p}$ and $\|\cdot\|_{2\to q}$ in~\cite{FS72, IS70}.}
Semenov and Shneiberg showed that if $T$ is any operator with $\|T\|_{p\to q} < \infty$
then for all $\epsilon < \epsilon_0 = \epsilon_0(p, q, \|T\|_{p\to q})$ we have
	$$ \|(1-\epsilon)\EE + \epsilon T\|_{p\to q} = 1\,,$$
provided that $\EE \circ T = T \circ\EE$, $T1 = 1$ and $(\EE f)(x) \eqdef \EE[f(X)]$. The key point is that $\epsilon_0$
only depends on $T$ through the norm $\|T\|_{p\to q}$. Paired with our
Theorem~\ref{th:main} this allows to establish that certain
permutation-invariant (or $S_n$-equivariant) operators in Hamming space have $L_p\to L_q$ norm equal to 1.

\subsection{Application: sumsets in Hamming space}
Our original interest in hypercontractivity was motivated by a remarkably simple solution it yields
to a problem that the author attempted to solve using more conventional semi-definite
 programming (SDP), compare Sections IV in~\cite{YP12-isit} and~\cite{YP13-htstruct}. Here is an application of the
new result (Theorem~\ref{th:main}) similar in spirit:
\begin{corollary} For every $\epsilon\in(0,1)$ there are constants $C_1, C_2 > 0$ such that for any
dimension $n$ and any set 
$A\subset \FF_2^n$ we have
	$$ \sup_{j\in[\epsilon n, (1-\epsilon)n]} {2^n(1_A * 1_A, 1_{\Sph_j})\over |\Sph_j| |A|} \ge
\lambda \quad \implies \quad |A| \ge C_1 \lambda^{C_2} 2^n\,. $$
In other words, $\PP[X+Y\in A] \ge \lambda$ implies $|A| \ge C_1 \lambda^{C_2} 2^n$, where $(X,Y)$ is uniform on $A\times
\Sph_j$.
\end{corollary}

\textbf{Remark:} It is known that any linear subspace $V \subset \FF_2^n$ which
contains a $\Omega(1)$-fraction of any $\Sph_{\delta n}$ must have co-dimension $O(1)$ (in
$n\to\infty$).\apxonly{\footnote{Here is how to prove this. First, we can show explicitly
that
 $$ \max_{\dim L = n-1} |L\cap \Sph_{\delta n}| = \bar \delta |\Sph_{\delta n}|$$ 
 This is shown by writing LHS as inner product and taking Fourier transform to get
 	$$ (1_{L^\perp}, K_{\delta n}) \to \max $$
	But since $1_{L^\perp} = \delta_0 + \delta_{v}$ we only need to optimize $v$ and
	the best choice is $|v|=1$, i.e. $L=\{x_1=0\}$. 
	Similarly for general co-dimension $k$ such that $2^k \approx \sum_0^r {n\choose
	j}$ we can show
	$$ \max_{\dim L = n-k} |L\cap \Sph_{\delta n}| \simleq |\Sph_{\delta n}|
	\EE[(1-2\delta)^W | W\le r]\,, \quad W \sim Bino(n, 1/2)\,.$$
	Here I upper-bounded $1_{L^\perp}$ by $1_{B_r}$ using monotonicity of Krawtchouks.
	Then I also approximated $K_{\delta n}(j) \simleq (1-2\delta)^j {n\choose \delta
	n}$.
	Overall, we see that for any fixed co-dimension $k$ the maximal fraction of
	$\Sph_{\delta n}$ decays as $k$ grows.}}%
{} This
corollary is a generalization: if a sumset $A+A$ contains a $\lambda$-fraction of
any Hamming sphere $\Sph_j$ (counted with multiplicity normalized by $|A|$) then the set must be of cardinality
$\Omega(2^n)$.

\begin{proof} We prove a stronger statement: 
\begin{equation}\label{eq:vnd}
	 \left(\phi*\phi, {1_{\Sph_j}\over |\Sph_j|}\right) \ge \lambda \|\phi\|_2^2 
	\quad \implies \quad {\|\phi\|_2^2\over
	\|\phi\|_1^2} \le {1\over C_1} \lambda^{-C_2}\,,
\end{equation}
from which the result follows by taking $\phi = 1_A$. To show~\eqref{eq:vnd} denote $\delta =
{j\over n}$ and consider the chain
\begin{align} \lambda \|\phi\|_2^2 &\le \left(\phi * \phi, {1_{\Sph_j}\over |\Sph_j|}\right)\\
			&= (\phi, S_\delta \phi)\\
			&\le \|\phi\|_2 \|S_\delta \phi\|_2\label{eq:la_1}\\
			&\le C \|\phi\|_2 \|\phi\|_p\,, \qquad p=1+(1-2\epsilon)^2 < 2
				\label{eq:la_2} \\
			&\le C \|\phi\|_2 \|\phi\|_1^{{2\over p} - 1} \|\phi\|_2^{2-{2\over p}}
				\label{eq:la_3} 
\end{align}
where~\eqref{eq:la_1} is Cauchy-Schwarz, ~\eqref{eq:la_2} is from Theorem~\ref{th:main},
and~\eqref{eq:la_3} is from log-convexity of ${1\over p} \mapsto \|\phi\|_p$. Rearranging terms
yields~\eqref{eq:vnd}.
\end{proof}

\def\FR{\mathtt{FR}}
In fact, this corollary can be interpreted in terms of the Frankl-R\"odl graphs $\FR_\gamma^n$, which are defined on the
vertex set $\FF_2^n$ with $v\sim v'$ if $|v-v'|=(1-\gamma)n$. Denoting by $E(A,A)$ the number of internal edges of a set
$A$, our corollary says
$$ |A| \le \mu 2^n \implies E(A,A) \le C_1' \mu^{C_2'} |\Sph_{\gamma n}| \, |A|\,.$$
In the regime of constant $\mu$ this is essentially tight. Indeed, an estimate in the opposite direction has been obtained by Benabbas, Hatami and Magen~\cite{benabbas2012isoperimetric}
(see~\cite[Section 5]{kauers2014hypercontractive} for a public account of these results):
\begin{equation}\label{eq:bhm_1}
	|A| \ge \mu 2^n \implies E(A,A) \ge \left((\mu/2)^{1\over \gamma} - o_n(1)\right) 2^n |\Sph_{\gamma n}|\,,
\end{equation}
provided $\gamma < 1/2$.
In particular, this implies that if $A$ is an independent set of $\FR_\gamma^n$ (so that $E(A,A)=0$) we must have $|A|
\le o(1) 2^n$. This is a weak form of the famous Frankl-R\"odl theorem~\cite{frankl1987forbidden} showing that
$\alpha(\FR_\gamma^n) \le (2-\epsilon(\gamma))^n$, where $\alpha(\cdot)$ denotes the maximal independent set of the
graph. Similar to our result,~\eqref{eq:bhm_1} was obtained by employing a reverse
hypercontractivity result of Borell~\cite{CB82}, which states
\begin{equation}\label{eq:revhc}
	\| N_\delta f\|_q \ge \|f\|_p, \qquad \forall -\infty<q<p<1, p-1 \le (q-1)(1-2\delta)^2\,, 
\end{equation}
for any $f>0$. Note that~\eqref{eq:revhc} cannot be extended to $S_\delta$, but in~\cite{benabbas2012isoperimetric} the
authors show that the eigenvalues of $N_\delta$ and ${1\over
2}(S_{\delta}+S_{\delta+1/n})$ are similar enough that the latter operator is almost reverse-hypercontractive. We will further discuss
results of~\cite{benabbas2012isoperimetric} below.

\subsection{Hypercontractivity and SDP} 

Part of our motivation to study hypercontractivity is that it may be employed as an 
improvement to the method of semi-definite programming (SDP) relaxation in various constraint satisfaction problems. 
For example, the best known bound~\cite{MRRW77} on the size of error
correcting codes in Hamming space are obtained by the SDP relaxation of Delsarte~\cite{PD73}, and 
there has long been interest in using hypercontractivity to improve the SDP relaxation, see~\cite{KL95}.

The relation between hypercontractivity and SDP has also been known in the computer science literature.\footnote{This
paper was originally written before some of the discussed results were published. We thank the reviewers for pointing
out these references.} For example,~\cite{GMPT10}
shows that any (fixed) number of rounds of Lov\'asz-Schrijver SDPs is unable to prove a bound better than
$\alpha(\FR(m,\gamma)) < ({1\over 2} - \epsilon) 2^m$, whereas we know from~\cite{frankl1987forbidden} that
$\alpha(FR(m,\gamma)) < (2-\epsilon)^m$. At the same time,~\cite{benabbas2012isoperimetric} shows that reverse
hypercontractivity proves $\alpha(FR(m,\gamma))< o(2^m)$. Following up on the latter,~\cite{kauers2014hypercontractive} 
shows that reverse hypercontractivity itself is provable in a sum-of-squares (SOS) proof system, thereby showing that
$\alpha(\FR(m,\gamma))<o(2^m)$ is provable via Lasserre's SOS algorithm of a fixed (but dependent on $\gamma$) degree.

This section gives another example where (direct, as
opposed to reverse) hypercontractivity supersedes SDP methods. We mention that while the previously mentioned 
examples deal with integer-programming problems, our example below is inherently ``continuous''.

Define $B_\delta(x) =
\delta^{|x|}(1-\delta)^{n-|x|}$ to be a distribution function of an iid Bernoulli noise. For $\lambda \in (0,1)$ we define
\begin{equation}\label{eq:la_0}
	 V_n(\lambda) = \max\left\{{(\phi, \phi)\over (\phi, 1)^2}: \phi \ge 0, 
				(\phi * \phi, B_\delta) \ge \lambda \|\phi\|_2^2 \right\}
\end{equation}
An argument entirely similar to~\eqref{eq:la_1}-\eqref{eq:la_3} invoking
Bonami-Gross~\eqref{eq:bonami} instead of Theorem~\ref{th:main} demonstrates\footnote{
The original question was to check whether there exists a small set $A\subset \FF_2^n$ such that  
$ \PP[X + X' = Z] \ge \lambda \PP[X+X' =0] $, where $X\dperp X' \sim$ uniform on $A$ and $Z\sim \Bern(\delta)$.
Bound~\eqref{eq:vnd2} shows any such set occupies a non-vanishing fraction of $\FF_2^n$.}
	\begin{equation}\label{eq:vnd2}
		V_n(\lambda) \le \lambda^{-s} 
\end{equation}	
for some $s>0$ and all dimensions $n$. 

Note that the problem in~\eqref{eq:la_0} is completely ``$L_2$'' and thus escaping to $L_p$
space in order to solve it looks somewhat unusual. Indeed, a more natural approach (at least to us) would be
to apply Fourier analysis or an SDP relaxation. Here is the ``spectral gap'' type of argument:
Since the second-largest eigenvalue of $N_\delta$ equals
$(1-2\delta)$ we get
	$$ (\phi_0, N_\delta \phi_0) \le (1-2\delta) \|\phi_0\|^2\,,$$
where $\phi_0 = \phi - (\phi, 1)$. Simple manipulations then imply
	$$ V_n(\lambda) \le {2\delta\over \lambda - (1-2\delta)}\,, \qquad \mathrm{if~}\lambda >
(1-2\delta)\,.$$
This proves a correct estimate of $O(1)$ but only for large values of $\lambda$.

An improvement of this method comes with the use of an SDP relaxation. The latter is obtained by considering $\psi = \phi*\phi$
and retaining only the non-negative definiteness property of $\psi$. I.e. we have the following upper bound:
	$$ V_n(\lambda) \le SDP(n, \lambda) \eqdef 
		\max \left\{2^n {(\psi, B_0)\over (\psi, 1)}: 
					\psi \ge 0\,, \psi \succeq 0\,, (\psi, B_\delta) \ge \lambda(\psi, B_0)\right\}
$$
where $B_0(x) = 1\{x=0\}$ and $\psi \succeq 0$ denotes that $f \mapsto f*\psi$ is a non-negative
definite operator. It can be shown that\footnote{These observations were made in collaboration with 
Prof. A. Megretski.}
$$ SDP(n, \lambda) = O(1)\,, \qquad \lambda > (1-2\delta)^2\,, $$
while for smaller values of $\lambda$ $SDP(n, \lambda)$ grows polynomially in $n$. Thus, while SDP improves on the
``spectral-gap'' argument, it is still unable to yield the correct estimate of $V_n(\lambda)$ for the entire range of $\lambda$.

\section{Auxiliary results}

\subsection{Notation}

For $x=(x_1,\ldots,x_n)\in \FF_2^n$ define $\bar x\eqdef(1-x_1,\ldots,1-x_n)$. 
For each $j=1,\ldots,n$ let
$$ \chi_j(x_1,\ldots,x_n) \eqdef 1_{\{x_j=0\}} - 1_{\{x_j=1\}}\,. $$
Define the characters, indexed by $v\in\FF_2^n$,
$$ \chi_v(x) \eqdef \prod_{j: v_j = 1} \chi_j(x) = (-1)^{\langle v,x \rangle}\,,$$
where $\langle v,x \rangle=\sum_{j=1}^n v_j x_j$ is a non-degenerate bi-linear form on $\FF_2^n$. The
Fourier transform of $f:\FF_2^n \to \mcomplex$ is 
$$ \hat f(\omega) \eqdef \sum_{x\in\FF_2^n} \chi_{\omega}(x) f(x) = 2^n (f, \chi_\omega)\,, \qquad
\omega \in \FF_2^n\,.$$

$L_p$-norms are monotonic
\begin{equation}\label{eq:lp_mono}
	 \|f\|_p \le \|f\|_{p_1}, \qquad p\le p_1\,.
\end{equation}
and satisfy the Young inequality:
\begin{equation}\label{eq:young}
	 \|f * g\|_p \le 2^n \|f\|_q \|g\|_r\, \qquad {1\over p} + 1 = {1\over q} + {1\over r}\,, 1\le
p,q,r \le \infty 
\end{equation}
\apxonly{Other useful identities:
\begin{align} (f, g) &= 2^{-n} (\hat f, \hat g) \\
	 	(f, \hat g) &= (\hat f, g) \\
		\widehat{ \widehat f} &= 2^n f \\
	   \|f\|_2 &= 2^{-{n\over 2}} \|\hat f\|_2\\
	   \|\hat f\|_{p'} &\le 2^{n\over p} \|f\|_p, \qquad 1 \le p \le 2\\
	   \widehat{B_\delta}(\omega) &= (1-2\delta)^{|\omega|}\,,\qquad B_\delta(x) =
\delta^{|x|}(1-\delta)^{n-|x|}\\
		\widehat{f*g} &= \hat f \cdot \hat g\\
		\widehat{fg} &= 2^{-n} \hat f * \hat g\\
		f'(x) &\eqdef \overline{f(-x)}\qquad\text{involutive anti-automorphism}\\
		(a, b*c) &= (a*c', b), \qquad c'(g) \eqdef \bar c(g^{-1})
\end{align}

 Relation to Walsh decomposition:
 	\begin{align} f(x) &= 2^{-n} \sum_\omega \hat f(\omega) (-1)^{\langle \omega, x \rangle} = \sum_{S\subset [n]} \hat f(S) W_S(x)\\
	W_S(x) &\eqdef W_S(x_S) = \prod_{i\in S} (-1)^{x_i},  \quad \hat f(S) \eqdef 2^{-n} \hat f(\omega), \omega \leftrightarrow S\\
	   \EE[f(X)|X_S=x] &= \sum_{S' \subseteq S} \hat f(S') W_{S'}(x)\\
	   L_i f(x) &\eqdef {1\over2}(f(x)-f(x+e_i)) = \sum_{S\supset\{i\}} \hat f(S) W_S(x)\\
	   L_i L_j f(x) &= \sum_{S\supset\{i,j\}} \hat f(S) W_S(x)	\qquad \mbox{($L_i$ commuting projectors!)}
\end{align}	   
}
For the size of Hamming spheres we have
\begin{equation}\label{eq:sphere}
	 |\Sph_{\delta n}| = {n \choose \lfloor \delta n\rfloor} = e^{n h(\delta) - {1\over 2} \ln n +
O(1)}\,, \qquad n\to\infty
\end{equation}
where the estimate is a consequence of Stirling's formula, $O(1)$ is uniform in $\delta$ on compact
subsets of $(0,1)$ and
\begin{equation}\label{eq:hbin}
	 h(\delta) = -\delta \ln \delta - (1-\delta) \ln (1-\delta)\,.
\end{equation}
Furthermore, for all $0 \le j \le n$
	\begin{equation}\label{eq:nck_1}
	e^{n h({j\over n})} \sqrt{1\over 2n} \le |\Sph_j| < e^{n h({j\over n})} 
\end{equation}
and for $1 \le j \le n-1$, cf.~\cite[Exc. 5.8]{RG68},
\begin{align}\label{eq:nck_2}
	e^{n h({j\over n})} \sqrt{n\over 8 j (n-j)} \le |\Sph_j| \le e^{n h({j\over n})}
		\sqrt{n\over2\pi j (n-j)}
\end{align}
\apxonly{Traditionally, the lower bound is given with $8$ replaced by $2\pi*e^{-1/3}$, but Gallager
tightens it by calculating for $j=1,2$ explicitly.

Other useful estimates:
\begin{align} 
	|\Sph_{\alpha n}| &= {\exp (n h(\alpha))  \over \sqrt{2 \pi n \alpha (1-\alpha)}} \cdot (1+O(n^{-1}))\,, \\
  	   | \mathbb{B}_{\alpha n}| & = |\Sph_{\alpha n}| \cdot {1-\alpha\over 1-2\alpha} \cdot (1+O(n^{-1}))\,, \qquad
	   \alpha < {1\over2}\,.
\end{align}	   
}

\subsection{Asymptotics of Krawtchouk polynomials}

Krawtchouk polynomials are defined as Fourier transforms of Hamming spheres:
\begin{equation}\label{eq:kr_sp}
	 K_j(x) \eqdef \widehat{1_{\Sph_j}} (x) = \sum_{k=0}^n (-1)^k {|x| \choose k} {n-|x| \choose
j-k}
\end{equation}
Since $K_j(x)$ only depends on $x$ through its Hamming weight $|x|$, we will abuse notation and
write $K_j(2)$ to mean value of $K_j$ at a point with weight $2$, etc.

Some useful properties of $K_j$, cf.~\cite{KL01}:
\begin{align}\label{eq:kr_sym}
	 K_j(x) &= (-1)^j K_j(n-x)\\
	 K_j(x) &= (-1)^x K_{n-j}(x)\label{eq:kr_jsym}\\
	 {K_j(x) \over K_j(0)} &= {K_x(j) \over K_x(0)}\label{eq:kr_exch}\\
	 K_j(0) &= \|K_j\|_2^2 = |\Sph_j| = {n \choose j}\,,\\
	 K_j(x) &= \sum_{|v|=j} \chi_v(x)
\end{align}

It is also well-known that $K_j(x)$ has $j$ simple real roots. For $j\le n/2$ all of them are in the
following interval, see~\cite[eq. (71)]{KL01}:
	$$ {n\over 2} - \sqrt{j (n-j)} \le x \le {n\over 2} + \sqrt{j (n-j)}\,.$$
For large $n$ the above bounds become tight, so that for $j=\delta n$ the location of the first root is at roughly
$$ \xi_{crit}(\delta) \eqdef {1\over2} - \sqrt{\delta(1-\delta)}\,.$$
\apxonly{And we have
\begin{align} \xi_{crit}(\delta) \le \xi \le 1-\xi_{crit}(\delta) \quad &\iff\quad 
	\xi_{crit}(\xi) \le \delta \le 1-\xi_{crit}(\xi) \\
	&\iff (1-2\xi)^2 + (1-2\delta)^2 \le 1 
\end{align}}

The following gives a convenient non-asymptotic estimate of the magnitude of $K_j(x)$:
\begin{lemma}\label{th:kasymp} For all $x,j=0,\ldots,n$ we have
\begin{equation}\label{eq:kasymp}
	 |K_j(x)| \le e^{n E_{j/n}(x/n)}\,,
\end{equation}
where the function $E_\delta(\xi) = E_{1-\delta}(\xi)$  and for $\delta\in[0,1/2]$:
\begin{equation}\label{eq:ed_2}
	 E_\delta(\xi) = \begin{cases}
			{1\over 2}\left(h(\delta) + \ln 2 - h(\xi)\right)\,,&
				\xi_{crit}(\delta) \le \xi \le 1-\xi_{crit}(\delta)\\
			\phi(\xi, \omega)\,, & \xi = {1\over2}( 1-(1-\delta)\omega - \delta
\omega^{-1})\,,  \end{cases}
\end{equation}
where in the second case $\omega$ ranges in
	$$ \omega \in \left[-\sqrt{\delta\over 1-\delta}, -{\delta\over 1-\delta}\right] \cup 
			\left[{\delta\over 1-\delta}, \sqrt{\delta\over 1-\delta}\right] $$
and
	\begin{equation}\label{eq:phidef}
	 \phi(\xi, \omega) \eqdef  \xi \ln|1-\omega| + (1-\xi) \ln|1+\omega| - \delta
\ln|\omega|\,.
\end{equation}
\end{lemma}

\textbf{Remark:} Exponent $E_\xi(\delta)$ was derived in~\cite{KL95} for $\xi\le\xi_{crit}(\delta)$.
Subsequently, a refined asymptotic expansion for all $\xi\in[0,1]$ was found in \cite{IS98}:
\begin{equation}\label{eq:kl_exact}
	 K_{\delta n}(\xi n) = {O(1)\over \sqrt{n}} e^{nE_\delta(\xi)}\,, 
\end{equation}
where the $O(1)$ term is $\theta(1)$ for $\xi \le \xi_{crit}$, while for
$\xi\in[\xi_{crit}, 1/2]$ the factor $O(1)$ is oscillating and may
reduce the exponent for a few integer points $x\in[\xi_{crit}n, (1-\xi_{crit})n]$,
which are close to one of the roots of $K_j(\cdot)$.
\begin{proof} Following~\cite{IS98}\footnote{Note that $K_j(\cdot)$
in~\cite{IS98} corresponds to $(-1)^j K_j(\cdot)$ in this paper.} we have
\begin{equation}\label{eq:ed_0x}
	 K_j(x) = {1\over 2\pi i} \oint_{\matc} (1-z)^x (1+z)^{n-x} z^{-j} {dz\over z}\,,
\end{equation}
where integration is over an arbitrary circle $\matc$ with center at $z=0$.
The derivative of the function $(1-z)^x (1+z)^{n-x} z^{-j}$ is zero when
\begin{equation}\label{eq:ed_0}
	 n-2x = (n-j) z + j z^{-1}\,. 
\end{equation}
Due to~\eqref{eq:kr_jsym} it is sufficient to consider $j\le n/2$.
Among the two solutions of~\eqref{eq:ed_0} denote by $\omega$ the unique one with smallest $|z|$ and $\Im(z)\ge 0$. Set,
for convenience
	$$ \xi = x/n, \qquad \delta = j/n \in[0,1/2] $$
and note that we have the following relation between $\omega$ and $\xi$
\begin{align}
	 \omega &= {1\over 2(1-\delta)} \left(1-2\xi - \sign(1-2\xi)\cdot \sqrt{(1-2\xi)^2 - 1 +
(1-2\delta)^2}\right) \\
	 1-2\xi &= (1-\delta)\omega + {\delta\over \omega}\,.\label{eq:xi_om}
\end{align}
	\apxonly{Other useful identities:
		\begin{enumerate}
			\item $$1-2\xi = 1- {(1-\omega)(\bar\delta \omega-\delta)\over \omega} $$
			\item Derivative of $\phi(\xi, \omega(\xi))$ under fixed $\delta$ is awesome:
				$$ {d\over d \xi} \phi(\xi, \omega(\xi)) =
				\log{1-\omega(\xi)\over1+\omega(\xi)}\,.$$
				Similarly for parametrization in terms of $\omega$
				$$ {d\over d\omega} \phi(\xi(\omega), \omega) =\xi'(\omega)  \log{1-\omega\over 1+\omega}
				= {1\over2}{\beta - (1-\beta)\omega^2\over \omega^2}\log{1-\omega\over
				1+\omega} $$
			\end{enumerate} }
As $\xi$ ranges from $0$ to $1$ the saddle point $\omega$ traverses the path 
	$$ \omega: {\delta\over 1-\delta} \to \sqrt{\delta \over 1-\delta} \to
-\sqrt{\delta\over1-\delta} \to -{\delta\over 1-\delta}\,,$$
where the middle segment is along the arc $e^{i \phi} \sqrt{\delta \over 1-\delta}, \phi\in[0,\pi]$;
Corresponding to these corner points $\xi$ ranges as follows
	$$ \xi: 0 \to \xi_{crit} \to 1-\xi_{crit} \to 1\,.$$
It is more convenient to reparameterize the answer in terms of location of the saddle-point
$\omega$. If we take $\matc$ to be the circle passing through $\omega$, then as
shown in~\cite[(3.4) and paragraph after (3.19)]{IS98} the maximum
	$$ \max_{z\in\matc} \left| (1-z)^x (1+z)^{n-x} z^{-j}\right|  $$
is attained at $z=\omega$ and is equal to $e^{n E_{\delta}(\xi)}$, where
\begin{equation}\label{eq:ed_1}
	 E_\delta(\xi) = \phi(\xi,\omega)\,,
\end{equation}
	and $\xi$ is a function of $\omega$ defined via~\eqref{eq:xi_om}. Thus, upper-bounding the
integrand $\{\cdot\}$ in~\eqref{eq:ed_0x} by the maximal value and noting that for any circle
	$$ \oint_\matc \left|{dz\over z}\right| \le 2\pi $$
we conclude that~\eqref{eq:kasymp} holds.

It remains to show the simplified expression in ~\eqref{eq:ed_2} for $\xi \in [\xi_{crit},
1-\xi_{crit}]$.
To that end, notice that such $\xi$ corresponds to 
	$$ \omega = e^{i \phi} \sqrt{\delta \over 1-\delta}, \qquad \phi \in [0,\pi]\,.$$
Substituting this $\omega$ into~\eqref{eq:ed_1} we see that~\eqref{eq:ed_2} is equivalent to
\begin{equation}\label{eq:ed_3}
	 \xi \ln {|1-\omega|\over \sqrt{\xi}} + (1-\xi) \ln {|1+\omega|\over \sqrt{1-\xi}} = {1\over
2} \ln {2\over 1-\delta}\,. 
\end{equation}
But for $\omega$ on the arc we have 
	$$ {|1-\omega|\over \sqrt{\xi}} = {|1+\omega|\over \sqrt{1-\xi}} = \sqrt{2\over
1-\delta}\,,$$
thus verifying~\eqref{eq:ed_3} and~\eqref{eq:ed_2}.
\end{proof}

\begin{figure}
\centering
\includegraphics[height=2in]{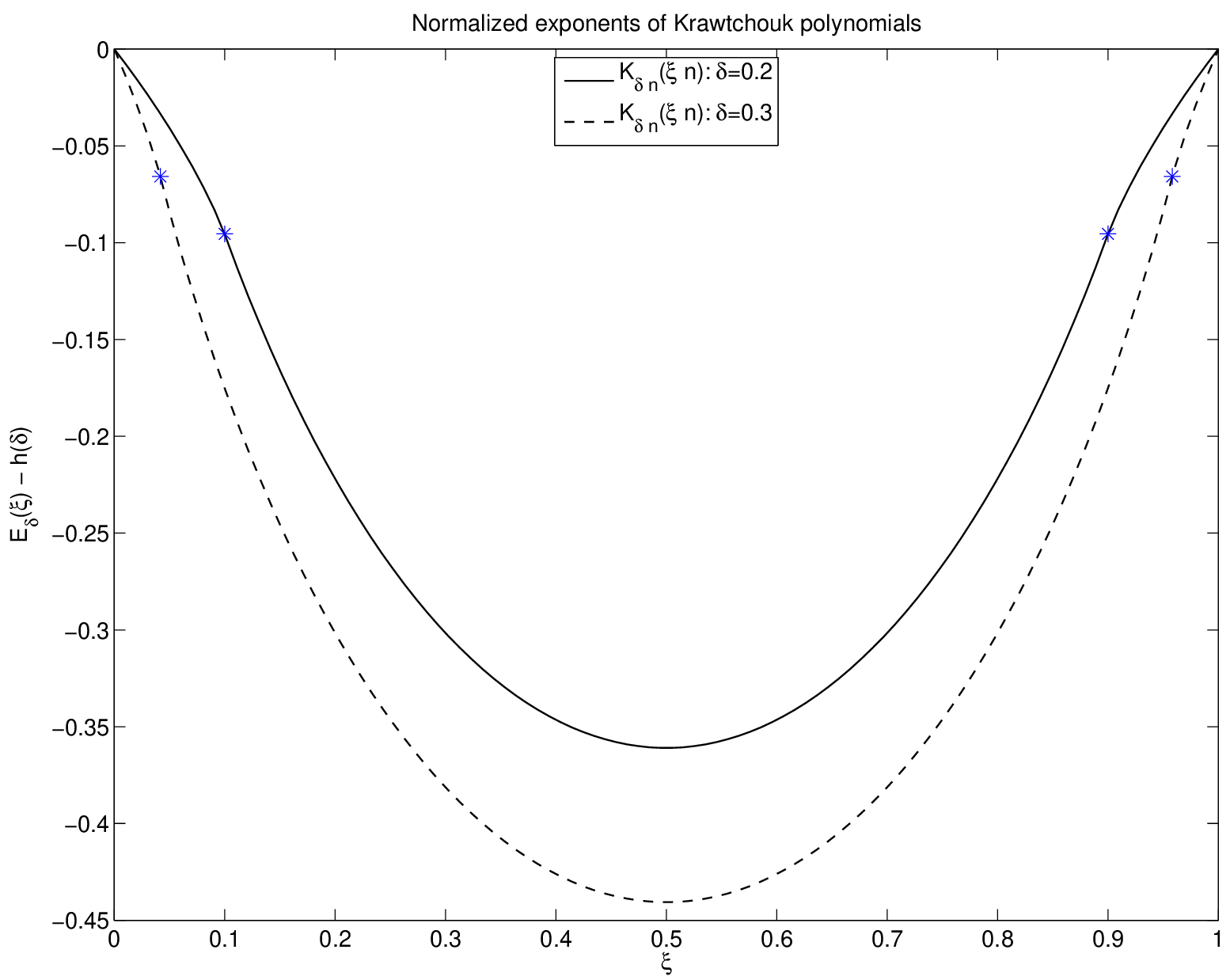}
\caption{The exponent of ${K_{\delta n}(\xi n)\over K_{\delta n}(0)}$ is equal
to $E_\delta(\xi)-h(\delta)$. The figure compares these exponents for two values of $\delta$.
Asterisks mark the interval $[\xi_{crit},1-\xi_{crit}]$ containing all the roots of $K_{\delta
n}(\cdot)$. In this interval $K_{\delta n}(\cdot)$ is oscillatory.}\label{fig:compare_2expon}.
\end{figure}

Some of the properties of $E_\delta(\xi)$ are summarized below (see Fig.~\ref{fig:compare_2expon} for an illustration):
\begin{enumerate}
	\item $(\delta,\xi) \mapsto E_\delta(\xi)$ is continuous on $[0,1]\times[0,1]$ and has
two symmetries: $E_\delta(\xi) = E_{1-\delta}(\xi)$, $E_\delta(\xi) = E_\delta(1-\xi)$.
	\item $E_\delta(0) = E_\delta(1) = h(\delta)$, $E_\delta(1/2) = h(\delta)/2$
	\item $E_{1/2}(\xi) = \ln 2 - h(\xi)/2$
	\item $E_\delta(\xi) = h(\delta) - h(\xi) + E_\xi(\delta)$
	\item $\xi \mapsto E_\delta(\xi)$ is monotonically decreasing on $[0,1/2]$ and 
		has continuous derivative on $[0,1]$.
	\item $\delta \mapsto E_\delta(\xi)$ is monotonically increasing on $[0,1/2]$.
	\item $\delta \mapsto E_\delta(\xi) - h(\delta)$ is monotonically decreasing on $[0,1/2]$.
	\item For fixed $\delta$ and all $\xi \le \xi_{crit}(\delta)$ we have
		\begin{equation}\label{eq:edx_a}
	 E_\delta(\xi) \le \xi \ln(1-2\delta) + h(\delta)\,.
\end{equation}
\end{enumerate}

We will also need a more refined estimate for $K_j(x)$ when $x$ is small:
\begin{lemma}\label{th:rdd} For $j\le n/2$ and $0\le x \le n\xi_{crit}(j/n) = n/2 - \sqrt{j(n-j)}$ we have
\begin{equation}\label{eq:rdd}
	 {K_j(x)\over K_j(0)} \le \left(1-{2j\over n}\right)^x\,.
\end{equation}
\end{lemma}
\textbf{Remark:} With the additional factor $O(\sqrt{n})$ the estimate~\eqref{eq:rdd} follows
from~\eqref{eq:kasymp}. Lemma~\ref{th:rdd} establishes the crucial relation 
between spectra of operators $N_\delta$ and $S_\delta$ powering Theorem~\ref{th:main}.
\begin{proof} In the mentioned range of $x$ the polynomial $K_j(x)$ is monotonically decreasing
since $K_j(0)>0$ and all roots are to the right of $x$. Hence, for any $x+1\le n\xi_{crit}(j/n)$ we have
\begin{equation}\label{eq:rdd0}
	 0 \le {K_j(x+1)\over K_j(x)} < 1\,.
\end{equation}
On the other hand, e.g.~\cite[(15)]{KL01}, $K_j(\cdot)$ satisfies a three-term recurrence
	\begin{equation}\label{eq:rdd_x}
	 (n-x) K_j(x+1) - (n-2j) K_j(x) + x K_j(x-1) = 0\,.
\end{equation}
Dividing by $nK_j(x)$ we get
\begin{align}
	 {K_j(x+1)\over K_j(x)} &= \left(1-{2j\over n}\right) - {x\over n}\left( {K_j(x-1)\over K_j(x)} -
{K_j(x+1)\over K_j(x)}\right)\label{eq:rdd1a} \\
	&\le \left(1-{2j\over n}\right)\,,\label{eq:rdd1}
\end{align}
where~\eqref{eq:rdd1} is from~\eqref{eq:rdd0}. The~\eqref{eq:rdd} then follows by
iterating~\eqref{eq:rdd1}.
\end{proof}

Note that for $j\approx {n\over 2}$ conditions of Lemma~\ref{th:rdd} are not satisfied for any $x$.
For such $j$ we prove another (somewhat loose) estimate below. 
\begin{lemma}\label{th:rqq} Fix arbitrary $\theta_1 \in (0,1/2)$. Then for all $x,j$ such that
\begin{align}
	 n-2j &\le n \theta_1,\label{eq:rqq_0}\\
	0\le x &\le 1+{\theta_1\over 1+\theta_1^2}(n\theta_1 - (n-2j))\label{eq:rqq_0a}
\end{align}
we have
\begin{equation}\label{eq:rqq}
	\left| {K_j(x) \over K_j(0)} \right| \le \theta_1^x\,.
\end{equation}
\end{lemma}
\begin{proof} Denote $\theta = 1-2{j\over n} \le \theta_1$. 
Clearly~\eqref{eq:rqq} holds for $x=0$. From~\eqref{eq:rdd1a} and~\eqref{eq:rqq_0} it also holds for
$x=1$. Let the induction hypothesis be that~\eqref{eq:rqq} holds for $x \le x_0$. Then 
	\begin{align} \left|K_j(x_0+1)\over K_j(0)\right| &= 
			\left| {n\theta\over n-x_0} {K_j(x_0)\over K_j(0)} - {x_0\over n-x_0}
{K_j(x_0-1) \over K_j(0)}\right| \label{eq:rqq_1}\\
			&\le {n \theta\over n-x_0} \theta_1^{x_0} + {x_0\over n-x_0}
\theta_1^{x_0-1}\,,\label{eq:rqq_2}
\end{align}
where~\eqref{eq:rqq_1} is from~\eqref{eq:rdd_x} and~\eqref{eq:rqq_2} is by induction hypothesis.
Finally, it is easy to see that whenever $n-x_0 > 0$ it holds that
$$ x_0 \le n {\theta_1 \over 1+\theta_1^2}(\theta_1 - \theta) \quad \iff \quad 
{n \theta\over n-x_0} \theta_1^{x_0} + {x_0\over n-x_0} \theta_1^{x_0-1} \le \theta_1^{x_0+1}\,,$$
which concludes the proof of~\eqref{eq:rqq} for $x=x_0+1$.
\end{proof}
\apxonly{
\textbf{Remark:} Selecting $\theta_1 > {\theta(1+\sqrt{5})\over 1+\sqrt{1-4\theta^2}}$ ensures
that~\eqref{eq:rqq_0a} includes some $x \ge n\xi_{crit}(j/n)$. In particular Lemma~\ref{th:rqq}
applies to a larger range of $x$ than Lemma~\ref{th:rdd} only when $n-2j \le {n\over \sqrt{5}}$.}

On the other extreme, for small values of $j$ we can extend Lemma~\ref{th:rdd} to the whole range $0
\le x\le {n\over 2}$:
\begin{lemma}\label{th:rzz} There exist $C_1\ge1$ and $\delta_0\in(0,1)$ such that for all $0 \le
j\le {\delta_0n}$ we have
	$$ \left|K_j(x)\over K_j(0)\right| \le C_1 \cdot \left(1-{2j\over n}\right)^x, \qquad 0 \le x \le
{n\over 2}\,.$$
\end{lemma}
\textbf{Remark:} In fact, one can show the statement with $C_1=1$ and $\delta_0=0.16$. This is
achieved by carefully following constants in the analysis and showing that
$\max_{\delta\in[0,\delta_0]}$ over the right-hand side of~\eqref{eq:nck_3} is $\le 1$ for $n\ge
300$.\apxonly{the rhs is $\le 1$ for $n\ge 156$.} For smaller $n$ 
the statement is checkable numerically, e.g. by running the recurrence~\eqref{eq:rdd_x} for normalized functions
${K_j(x)\over K_j(0)\left(1-2\delta\right)^x}$ (to avoid large numbers).
\begin{proof}For $j=0$ the inequality is trivial. For $x \le \xi_{crit}(j/n)$ it follows from Lemma~\ref{th:rdd}. Thus,
it is sufficient to consider $x \ge \xi_{crit}(j/n)$, $j\ge 1$. Denote $\delta = j/n$. Then from
Lemma~\ref{th:kasymp} and~\eqref{eq:nck_2} we have for all $n\ge 1$:
	\begin{equation}\label{eq:nck_3}
		\left|K_j(x)\over K_j(0)\left(1-2\delta\right)^x\right| \le
		\sqrt{8(1-\delta)} \cdot e^{n (f(\delta) - {1\over 2} h(\delta))} \sqrt{n \delta}\,,
\end{equation}		
	where
	$$ f(\delta) = \max_{\xi\in[\xi_{crit}(\delta), 1/2]} {1\over 2}(\ln 2 - h(\xi)) - \xi
\ln(1-2\delta)\,.$$
	From convexity of the function under maximization, we conclude 
	$$ f(\delta) = {\ln 2\over 2} - {1\over 2}\min \left( h(\xi_{crit}(\delta)) +
2\xi_{crit}(\delta) \ln (1-2\delta), \ln 2(1-2\delta)\right)\,.$$
	Taking derivative at $\delta=0$ we conclude that for some $\delta_0'>0$ we have
	$$ h(\xi_{crit}(\delta)) + 2\xi_{crit}(\delta) \ln (1-2\delta) \le \ln
		2(1-2\delta)\,,\qquad \forall\delta\in[0,\delta_0']\,.$$
	\apxonly{(Numerically: $\delta_0' \approx 0.1741$.)}
	Consequently, for such $\delta$
	$$ f(\delta) = {1\over 2}\left(\ln 2 - h(\xi_{crit}(\delta))\right) - \xi_{crit}(\delta)
\ln(1-2\delta)\,.$$
	Evidently, $f$ is continuously differentiable and
		$$ f(\delta) = 2\delta + o(\delta), \qquad \delta\to0\,.$$
	Therefore for some $\delta_0\in(0,\delta_0']$ we must have
		$$ f(\delta) - {1\over 2} h(\delta) < 0\,, \qquad \forall \delta\in(0,
\delta_0]\,.$$
	\apxonly{(Numerically: $\delta_0 = \delta_0'$.)}
	The statement of the Lemma then follows with $C_1 = \max(1, \sqrt{8}C_1')$,
where $C_1'$ is the finite supremum found in the following Lemma.
\end{proof}

\begin{lemma} Let $\alpha,\delta_0,C>0$
and $f$ -- a continuous function on $[0, \delta_0]$ with $f(0)=0$, derivative (one-sided at
$0$) bounded by $C$ and satisfying
\begin{equation}\label{eq:rzz_1}
	 f(\delta) - \alpha h(\delta) < 0\,, \qquad \forall \delta \in(0,\delta_0]\,.
\end{equation}
Then 
	\begin{equation}\label{eq:rzz}
		\sup_{n\ge 1} \max_{\delta \in[0, \delta_0]} e^{n (f(\delta) - \alpha h(\delta))}
\sqrt{n\delta} < \infty\,.
\end{equation}
\end{lemma}
\apxonly{Remark: The bound is order-tight and extendable to (also tight)
	$$ \max_{n\ge 1, \delta \in[0, \delta_0]} e^{n (f(\delta) - \alpha h(\delta))}
\sqrt{n^\beta\delta} = \theta(n^{\beta-1})\,.$$
}
\begin{proof} Under conditions of the theorem there exists $0<\delta_1 < \delta_0$ such that
	$$ f(\delta) \le {\alpha\over 2} h(\delta)\,,\qquad \forall \delta\in[0,\delta_1]\,.$$
Thus we have
\begin{align} \max_{\delta\in[0,\delta_1]} n(f(\delta) - \alpha h(\delta)) + {1\over 2} \ln (\delta n) &\le
		{1\over 2} \max_{\delta\in[0,\delta_1]} -\alpha n h(\delta) + \ln \delta n\\
		&\le{1\over 2} \max_{\delta\in[0,\delta_1]} \alpha n \delta \ln \delta + \ln (\delta n)\,.
		\label{eq:rzz_2}
\end{align}
Without loss of generality we may assume $\delta_1 < {1\over e}$ and $n>{e^2\over \alpha}$. In this
case, maximization in~\eqref{eq:rzz_2} is attained at $\delta^* \in (0, {1\over n\alpha})$. 
\apxonly{To show this, first take derivative in $\delta$ and show that its minimum is attained at
$\delta_{flex} = {1\over n\alpha}$ and that this derivative is still negative at $\delta_1<{1\over
e}$. Hence, the zero-crossing is in the mentioned interval.}
Consequently, upper-bounding the first term by zero and second by $\ln({1\over n\alpha}\cdot n)$ we get
$$
	{1\over 2} \max_{\delta\in[0,\delta_1]} \alpha n \delta \ln \delta + \ln (\delta n) \le
	{-\ln \alpha\over 2}\,.$$
On the other hand, from~\eqref{eq:rzz_1} and continuity we get
	$$ \max_{\delta\in[\delta_1,\delta_0]} f(\delta) - \alpha h(\delta) = -C_2 < 0\,.$$
Therefore, putting both bounds together
	$$ \max_{n\ge 1, \delta \in[0, \delta_0]} e^{n (f(\delta) - \alpha h(\delta))} \le
		\max\left({1\over \sqrt{\alpha}}, \sup_n \sqrt{\delta_0 n} e^{-C_2 n}\right) <
\infty\,.$$
\end{proof}

\textbf{Remark:} Reference~\cite{benabbas2012isoperimetric} establishes the following estimate:
$$ \left| {1\over 2} \left( {K_c(n)\over K_c(0)} + {K_{c-1}(n)\over K_{c-1}(0)} \right) - \left(1-{2c\over
n}\right)^n\right| \le O(\max(n^{-{1\over 5}}, {n\over c^2} \log^2 {c^2\over n}))\,,$$
for all $e^2 \sqrt{n} c \le {n\over 2}$. This result is incomparable to ours: it bounds deviation from $(1-{2c\over
n})^n$ on both sides, albeit much less precisely.

Finally, for illustrating tightness of the bounds in the next section we will need the following Lemma, proved in the
Appendix. It is not used in the proof of Theorem~\ref{th:main}.
\begin{lemma}\label{th:krawt_lp} $L_p$ norms of Krawtchouk polynomials are given asymptotically by the following
parametric formula: Let $\omega \in [0,1]$ then for $p\ge 2$
	\begin{align} 
	\| K_{\lfloor\delta n\rfloor}\|_p &= \exp\left\{ n\left( {h(\xi) - \ln 2\over p} +
\phi(\xi,\omega)\right) + O(\log n)\right\}, \qquad n\to\infty\label{eq:klp}\\
		c &= {(1+\omega)^p -(1-\omega)^p\over (1+\omega)^p + (1-\omega)^p}\label{eq:klp_a}\\
		\xi &=  {1-c\over 2} = {1\over2}( 1-(1-\delta)\omega - \delta \omega^{-1})\\
		\delta &= {c\omega - \omega^2\over 1-\omega^2}\label{eq:klp_b}
\end{align}
	and $\phi(\xi,\omega)$ is given by~\eqref{eq:phidef}. For $p\le 2$ we have
	\begin{equation}\label{eq:klp_2}
		\| K_{\lfloor\delta n\rfloor}\|_p = \exp\left\{ {n\over 2} h(\delta) + O(\log n)
\right\}\,, 
\end{equation}
as $n\to\infty$ along a subsequence such that both $\lfloor \delta n\rfloor$ and $n$ are even.
\end{lemma}

\subsection{Norms of Fourier projection operators}
The Fourier projection operators $\Pi_a$ are defined as
\begin{equation}\label{eq:fproj}
	 \widehat{\Pi_a f} \eqdef \hat f \cdot 1_{\Sph_a}\,\qquad a=0,1,\ldots,n\,,
\end{equation}
or, equivalently,
	$$ \Pi_a f \eqdef 2^{-n} f*K_a\,.$$
On the other hand from Young's inequality~\eqref{eq:young} we have for any convolution operator:
	$$ \|\phi*(\cdot)\|_{1\to 2} = 2^n \|\phi\|_2 \,.$$
\apxonly{(For lower-bound $f=\delta_0$.)} Thus we have
	\begin{equation}\label{eq:fp_l1_exact}
	 \|\Pi_a\|_{1\to 2} = \sqrt{n \choose a}\,.
\end{equation}
Also, we note that
	$$ \|\Pi_a\|_{p\to q} = \|\Pi_{n-a}\|_{p\to q}\,,$$
and thus we only consider $a\le {n\over 2}$ below. \apxonly{Proof: $\widehat{f\cdot(-1)^{|x|}} = \hat
f(\bar \omega).$}

Estimates for other $L_p \to L_2$ follow from Bonami-Gross inequality~\eqref{eq:bonami} and complex
interpolation:
\begin{lemma}\label{th:pi_ub} For any $1 \le p \le 2$ and $0 \le a = n\delta \le {n\over2}$ we have 
	\begin{align} \|\Pi_a\|_{p\to 2} &\le \begin{cases}
						(p-1)^{-{a\over 2}}, &\qquad p > p^*\,,\\
						(p^*-1)^{-{(1-s)a\over 2}} {n\choose a}^{{s\over p} -
{s\over 2}}, & {1\over p} = {1-s\over p^*} + s, 0 \le s\le 1
					   \end{cases} 
		\label{eq:fp_best}
\end{align}
	where $p^* = p^*(a) = 2$ if ${h(\delta)\over \delta} \le 2$\apxonly{\footnote{Base
independent version: 
	$h(\delta) \ge 2 \delta \log e$. Threshold value: $\delta \approx 0.30925$.}}, and otherwise
$p^* \in(1,2)$ is a solution of
	$$ p^* - \ln (p^*-1) = {h(\delta) \over \delta}\,. $$
	We also have two weaker bounds
	\begin{align} \|\Pi_a\|_{p\to 2}  & \le (p-1)^{-{a\over 2}}\,,\label{eq:fp_bonami}\\
			\|\Pi_a\|_{p\to 2}&\le {n \choose a}^{{1\over p} - {1\over2}}\,. \label{eq:fp_l1}
	\end{align}
\end{lemma}
\textbf{Remark:} The estimate~\eqref{eq:fp_bonami} has been the basis of Kahn-Kalai-Linial
results~\cite{KKL88}, so we refer to~\eqref{eq:fp_bonami} as KKL bound. Note that $p^*(a) = 2$
corresponds to $a > 0.3093 n$, and then bound~\eqref{eq:fp_best} coincides with~\eqref{eq:fp_l1}.
\begin{proof} From Riesz-Thorin interpolation~\cite[Section VI.10.8]{DS58}, we know that the map ${1\over p} \mapsto \|\Pi_a\|_{p\to 2}$ is log-convex. Thus~\eqref{eq:fp_best}
follows from~\eqref{eq:fp_bonami} and~\eqref{eq:fp_l1} by convexification (the value of $p^*$ is
chosen to minimize the resulting exponent when $a=\delta n$). Thus, it is sufficient to
prove~\eqref{eq:fp_bonami} and~\eqref{eq:fp_l1}. The second one again follows from interpolating
between~\eqref{eq:fp_l1_exact} and $\|\Pi_a\|_{2\to 2} = 1$. For the first one notice that for any
$\tau$ we have
	$$ N_\tau \Pi_a = \Pi_a N_\tau = (1-2\tau)^a \Pi_a\,.$$
	And thus from~\eqref{eq:bonami} with $(1-2\tau)^2 = p-1$ we get
	$$ \|\Pi_a f\|_2 = |1-2\tau|^{-a} \|\Pi_a N_\tau f\|_2 \le |1-2\tau|^{-a} \|N_\tau f\|_2 \le
|1-2\tau|^{-a} \|f\|_p\,.$$
\end{proof}

\begin{figure}[t]
\centering
\includegraphics[width=.4\textwidth]{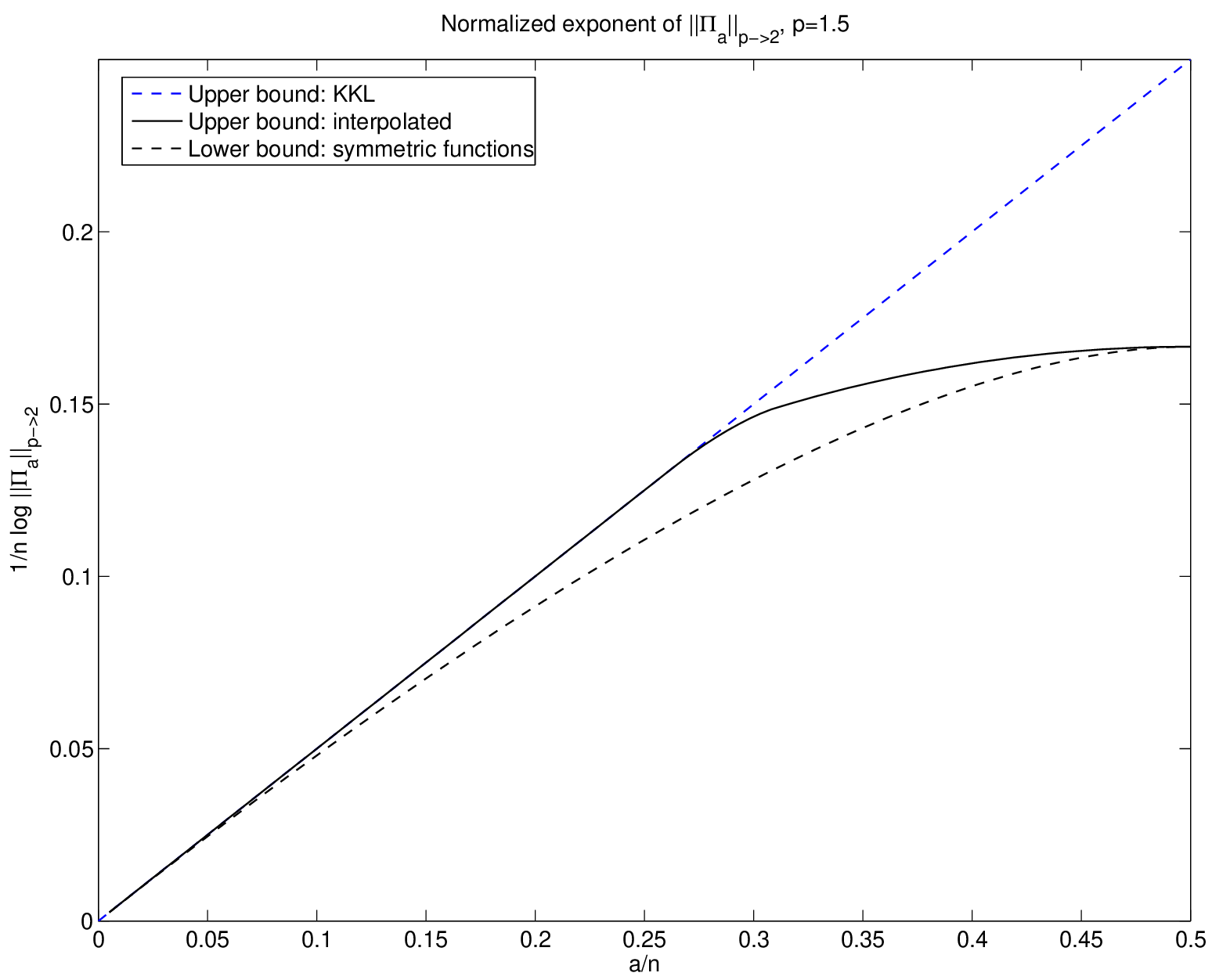}\hskip30pt
\includegraphics[width=.4\textwidth]{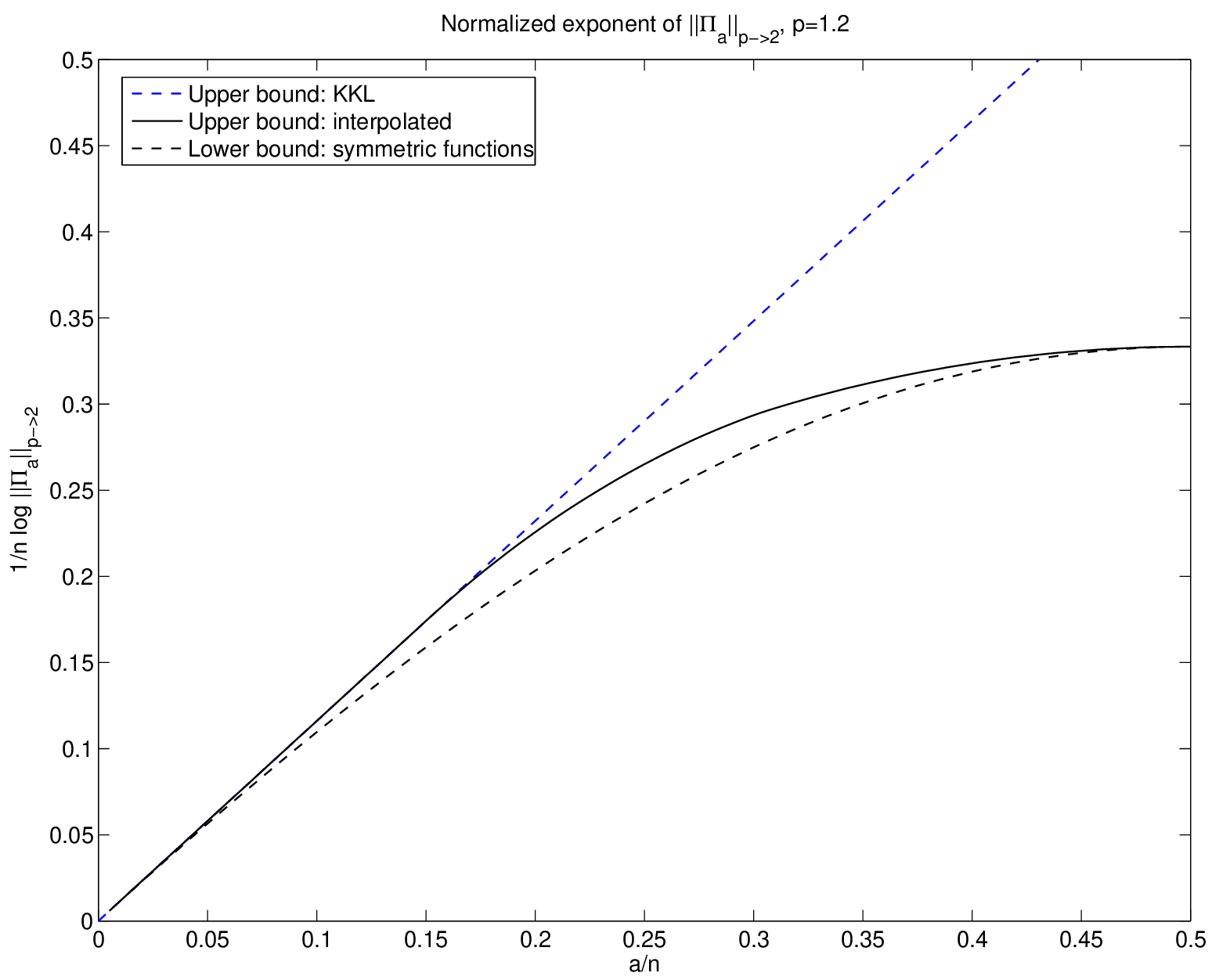}
\caption{Exponent of $\|\Pi_a\|_{p\to 2}$ as a function of $a$ for two values of $p$. Two upper
bounds correspond to Kahn-Kalai-Linial~\eqref{eq:fp_bonami} and the interpolated
one~\eqref{eq:fp_best}. The lower bound is given by
considering only permutation invariant functions (cf. Lemmas~\ref{th:krawt_lp} and~\ref{th:pi_lb}).}
\label{fig:compare_pa}.
\end{figure}

To verify the tightness of our bounds we derive a simple lower bound by considering permutation
invariant functions:
\begin{lemma}\label{th:pi_lb} For any $a\in\{0,\ldots,n\}$ and any $q,p \ge 1$ we have
	$$ \|\Pi_a\|_{p\to q} \ge {\|K_a\|_q \|K_a\|_{p'}\over \|K_a\|_2^2}\,,$$
	where $p'={p\over p-1}$ is the H\"older conjugate.
\end{lemma}
\apxonly{\\\textbf{TODO:} Amazingly, this exponentially coincides with the iid lower bound, obtained by taking $f=\prod
f_1(x_i)$.}
\begin{proof} The lower bound is shown by optimizing over a class of permutation invariant functions
$$ f(x) = K_a(x) + \sum_{j \neq a}^n c_j K_j(x) \eqdef K_a(x) + \Phi(x)\,, $$
where $\Phi \perp K_a$. Note that 
	\begin{align} \inf_{\Phi \perp K_a} \|f\|_p &= \inf_{\Phi \perp K_a} \sup_{g: \|g\|_{p'} \le 1} 
		(K_a + \Phi, g) \label{eq:pl_1}\\
		&=\inf_{\Phi \perp K_a} \sup_{g-\text{sym.}: \|g\|_{p'} \le 1} (K_a + \Phi, g)
\label{eq:pl_1a}\\
		&= \sup_{g-\text{sym.}: \|g\|_{p'} \le 1} \inf_{\Phi \perp K_a} (K_a + \Phi, g)
				\label{eq:pl_2} \\
		&= \left(K_a, {K_a\over \|K_a\|_{p'}}\right) = {\|K_a\|_2^2 \over \|K_a\|_{p'} }
				\label{eq:pl_3}\,, 
\end{align}
where~\eqref{eq:pl_1} is by duality $(L_p)^* = L_{p'}$,~\eqref{eq:pl_1a} states the obvious fact
that supremization can be restricted to permutation-symmetric $g$,~\eqref{eq:pl_2} is by Kneser's minimax 
theorem~\cite{kneser1952theoreme} (for bi-affine function over $X\times Y$ with $X$ convex-compact, $Y$ convex and $f$
upper semi-continuous on $X$) and~\eqref{eq:pl_3} is because the inner $\inf$ can only be finite if $g$ belongs to the
one-dimensional subspace spanned by $K_a$, i.e. $g=c K_a$ for a suitable $c$.

Since $\Pi_a (K_a + \Phi) = K_a$ we conclude that
	$$ \|\Pi_a\|_{p\to q} \ge {\|K_a\|_q\over \inf_{\Phi\perp K_a} \|K_a + \Phi\|_p}
			= {\|K_a\|_q \|K_a\|_{p'}\over \|K_a\|_2^2} $$
as claimed.
\end{proof}

On Fig.~\ref{fig:compare_pa} we compare the upper and lower bounds on $\|\Pi_a\|_{p\to 2}$ as $a$ ranges from 0
to $n/2$ for two values of $p$. We note that KKL bound~\eqref{eq:fp_bonami} is significantly
suboptimal for small $p$ and large $a$. For example, for $a>0.3093n$ the bound~\eqref{eq:fp_l1}
is strictly better than KKL.

\apxonly{
\medskip
\textbf{Notes on norms of Fourier projection:}
\begin{itemize}
	\item Estimates of $\|\Pi_a\|_{2\to q}$ are related to Khintchine-like inequalities that state
	$$ \sup_{f} {\|\Pi_a f\|_q \over \|\Pi_a f\|_2} \le C\,.$$
In particular, for $a=1$ and $q\ge 3$ symmetric $f$ indeed achieves the maximum, cf.~\cite{UH82}. At
the same time for $2 < q < 3$ this is not true: for example, 
	$$ {1\over\sqrt{2}} \|\chi_1 + \chi_2\|_{2.5} > {1\over\sqrt{3}} \|\chi_1 + \chi_2 + \chi_3\|_{2.5}
	$$
	\textbf{Later:}~\cite{UB82} seems to prove stuff for $k=1,n\to\infty$. Some other paper (dont remember) proved
	that $(a_1,\ldots,a_n) \to \|\sum a_i B_i\|_p^p/(\sum a_i^2)^{p/2}$ is Schur-concave for $p\ge 3$. That's the
	result I meant to cite.

\item Kirshner-Samorodnitsky'2018 closed the gap for $\|\Pi\|_{2\to 4} = \|\Pi\|_{4/3 \to 2}$. 
They proved that with exponential precision
	$$ {\|\Pi_{x n}\|_{4} \over \|\Pi_{x n}\|_2} = e^{o(n)} {\|1_{\Sph_{xn}}\|_4\over \|1_{\Sph_{xn}}\|_2}\,.$$
	This implies the $\|\Pi\|_{2\to4}$ via $\|\Pi f\|_4 \le e^{nE} \|\Pi f\|_2 \le e^{nE} \|f\|_2$.

\item Bourgain showed that in general for all $1\le p \le \infty$ (I
	doubt $p=1,\infty$ is actually true) and for all spaces we have
		$$ \|\Pi_k f\|_p \le c_p^k \|f\|_p\,, $$
		where $\Pi_k$ is orthogonal projection on the eigenspace corresponding to
		$\rho^k$ of $(\rho I + (1-\rho) \EE_\pi)^{\otimes n}$. For
		$\pi(x)=\pi(-x)$ this follows from KKL and $\pi=Bern(1/2)$ case. For
		non-symmetric $\pi$ there is the amazing symmetrization idea:
			$$ \tilde f(b^n, x^n) = \sum_{S\subset[n]} b_S f^{=S}(x_S) $$

\item Chang's Lemma states that for any subset $S$ with $\PP[S]=\alpha$ we have
	$$ \|\Pi_1 1_S\|_2^2 \le \alpha^2 \log {1\over \alpha^2}\,.$$
	This follows from~\eqref{eq:fp_bonami} with $p=1+{1\over
	\log {1\over\alpha}}$.
\item Lemma~\ref{th:pi_lb} and~\eqref{eq:klp_2} imply that for all $1\le p \le 2$ and $\delta > 0$ we have
	$$ e^{O(\log n)} \le \|\Pi_{\delta n}\|_{2\to p} = \|\Pi_{\delta n}\|_{p'\to 2} \le 1 \,.$$
	Since $\|\Pi f\|_2 \le \|f\|_2 \le \|f\|_{p'}$. Thus, $\Pi_{\delta n}$ are not exponentially-contracting in this regime.
\end{itemize}
}

Before proceeding to the proof of the main result, we need one last estimate relating magnitude of Krawtchouk
polynomials (in the oscillating strip) to the norms of projectors $\Pi_a$.

\begin{lemma}\label{th:tough} Fix arbitrary $0<\delta_0 < \Delta < 1/2$. Then there exist constants $C_1', C_2 > 0$ such that
for all $n\ge 1$, all $j\in[\delta_0 n, \Delta n]$ and all 
		$$ {n\over 2} - \sqrt{j(n-j)} \le x \le {n\over 2} + \sqrt{j(n-j)} $$
	we have
\begin{equation}\label{eq:tough}
	 \left|K_j(x) \over K_j(0)\right| \cdot \|\Pi_j\|_{p({j\over n})\to 2} \le C_1' \sqrt{n} e^{-C_2 n} 
\end{equation}
	where $p(\delta) = 1 + (1-2\delta)^2$.
\end{lemma}
\begin{proof} Let $\xi = {a\over n}$ and $\delta = {j\over n}$. From symmetry, we can and will assume 
$\xi \le {1\over 2}$. Since
$\xi$ is restricted to critical strip of Krawtchouk polynomial $K_{\delta n}$ from
Lemma~\ref{th:kasymp}, bound~\eqref{eq:nck_2} and Lemma~\ref{th:pi_ub} it is sufficient to show
\begin{equation}\label{eq:tg_1}
	\max_{\delta_0 \le \delta \le \Delta} \max_{\xi: (1-2\xi)^2 + (1-2\delta)^2 \le 1} {1\over2}
(\ln 2 -h(\xi) - h(\delta)) + \pi(p(\delta), \xi) \le -C_2 < 0\,,
\end{equation}
where $p(\delta) = 1+(1-2\delta)^2$ and
$$ {1\over p} \mapsto \pi(p, \xi) $$
is the convexification of the function (cf. Lemma~\ref{th:pi_ub})
\begin{equation}\label{eq:tg_0x}
	 {1\over p} \mapsto \min\left\{ -{\xi\over 2} \ln(p-1), ({1\over p} - {1\over 2})
h(\xi)\right\}\,.
\end{equation}
To show~\eqref{eq:tg_1} we first change variable $\delta$ to $p=p(\delta)=1+(1-2\delta)^2$. Set
	\begin{align} p_0 &= 1+(1-2\Delta)^2\,,\\
	   p_1 &= 1+(1-2\delta_0)^2\,. 
\end{align}
Then~\eqref{eq:tg_1} is equivalent to (we also interchange the maxima in $\xi$ and $\delta$):
\begin{equation}\label{eq:tg_2}
	\max_{\xi: (1-2\xi)^2 \le 2-p_0} \max_{p: p_0 \le p \le \min(p_1, 2-(1-2\xi)^2)} \eta(\xi,
p) + {\ln 2 - h(\xi)\over 2} \le -C_2 < 0
\end{equation}
where 
	$$ \eta(\xi, p) \eqdef \pi(p, \xi) - {1\over 2} h\left({1-\sqrt{p-1}\over 2}\right)\,.$$
By construction, ${1\over p} \mapsto \pi(p, \xi)$ is convex. Taking derivatives one can show that 
$h\left({1-\sqrt{p-1}\over 2}\right)$ is concave in $1\over p$. Thus, the maximization over $p$ in~\eqref{eq:tg_2}
is applied to a convex function and therefore must be achieved at one of the boundaries.
Consequently, to verify~\eqref{eq:tg_2} it is sufficient to show the following three strict
inequalities :
\begin{align} 
	\max_{\xi: (1-2\xi)^2 \le 2-p_0} \eta(\xi, p_0) + {\ln 2 - h(\xi)\over 2} &< 0
\label{eq:tg_3a}\\
	\max_{\xi: (1-2\xi)^2 \le 2-p_1} \eta(\xi, p_1) + {\ln 2 - h(\xi)\over 2} &< 0
\label{eq:tg_3b}\\
	\max_{\xi: 2-p_1 \le (1-2\xi)^2 \le 2-p_0} \eta(\xi, 2-(1-2\xi)^2) + {\ln 2 - h(\xi)\over 2}
&< 0 \label{eq:tg_3c} 
\end{align}
(the maximum value of the three left-hand sides is then taken to be $-C_2$).
The first two are verified as follows: From~\eqref{eq:tg_0x} we have
	$$ \pi(p,\xi) \le -{\xi\over 2} \ln(p-1)\,.$$
Plugging this upper bound in~\eqref{eq:tg_3a} we arrive at the optimization
	$$ \max_{\xi: (1-2\xi)^2 \le 2-p} -{\xi\over 2} \ln(p-1) - {1\over 2} h(\xi)\,.$$
Equating derivative in $\xi$ to zero, we find solution $\xi^*(p)=1-{1\over p}$. Since for $p>1$ we have $(1-2\xi^*(p))^2
< 2-p$ this is also the maximizer. Consequently, substituting $\xi=\xi^*(p)$ we get 
	$$ \max_{(1-2\xi)^2 \le 2-p} \eta(\xi, p) + {\ln 2 - h(\xi)\over 2} \le -{\xi^*(p)\over 2}
\ln(p-1)  + {1\over 2} \left[ \ln 2 - h(\xi^*(p)) - h\left({1-\sqrt{p-1}\over
2}\right)\right] $$
Function of a single variable $p$ on the right is continuous, non-positive 
and attains zero only at the endpoints of $p\in[1,2]$. Since both $p_0$ and $p_1$ belong to the interior of
$[1,2]$, this completes the proof of~\eqref{eq:tg_3a} and~\eqref{eq:tg_3b}.

To show~\eqref{eq:tg_3c} we apply the bound in~\eqref{eq:tg_0x} (without convexification):
\begin{equation}\label{eq:tg_5}
	 \max_{\xi} \eta(\xi, 2-(1-2\xi)^2) + {\ln 2 - h(\xi)\over 2} 
	\le \max_\xi {1\over 2} f(\xi) 
\end{equation}
where maximization is over 
\begin{equation}\label{eq:tg_6}
	 2-p_1 \le (1-2\xi)^2 \le 2-p_0  
\end{equation}
and $f(\xi)$ is defined as
\begin{align}
	 f(\xi) &\eqdef \min \left\{\left((1-2\xi)^2\over 2-(1-2\xi)^2\right)h(\xi), -\xi \ln
(4\xi(1-\xi))\right\} \nonumber\\
	&{} +  
		\ln 2 - h(\xi) - h\left({1\over 2} - \sqrt{\xi(1-\xi)}\right) 
\label{eq:tg_4}
\end{align}
The minimum in this expression selects the first term for $\xi \in[ \xi^*, 1/2]$ 
and second term otherwise, where $\xi^* \approx 0.3082$ is the solution of
$$ 8 \xi^2 (1-\xi) \ln \xi + (2\xi - (1-2\xi)^2) \ln (1-\xi) + 2\xi (2-(1-2\xi)^2) \ln 2 =0 $$
in the interior of $(0,1/2)$.
Furthermore, function in~\eqref{eq:tg_4} is
non-positive, continuous and attains zero only at $\xi=0,{1\over 2}$ both of which are excluded by the
constraints~\eqref{eq:tg_6}. Thus~\eqref{eq:tg_3c} holds.
\end{proof}

\section{Proof of Theorem~\ref{th:main}}\label{sec:proof}

Denote the boundary of $F$ as
	$$ p(\delta) \eqdef 1 + (1-2\delta)^2\,.$$
Note that every compact subset $K'$ of $F$ is contained in $F \cap \{p\ge p_0\}$ for sufficiently
small $p_0$ and in turn in some
\begin{equation}\label{eq:p_0}
	 K = \left(F \cap \{\delta: |1-2\delta| \ge \theta\}\right) \cup 
		\{(\delta,p):  |1-2\delta| \le \theta, p \ge p_0\} 
\end{equation}
for sufficiently small $\theta$. In particular, we may choose $\theta$ so small that
$p_0>1+\theta^2$. 
Next note that 
	$$ (f * 1_{\Sph_{n-a}})(x) = (f*1_{\Sph_a})(\bar x) $$
and thus estimates for $S_\delta$ and $S_{1-\delta}$ coincide asymptotically. 
Due to this symmetry and thanks to the monotonicity~\eqref{eq:lp_mono} of norms, to prove the theorem 
it is sufficient to prove the following pair of statements, corresponding to the boundary of $K$:
\begin{enumerate}
\item[S1.](critical estimate for $\delta<1/2$) For each $\delta$ there is $C_\delta$ such that for
all $n\ge 1$ and all functions $f$ we have
	\begin{equation}\label{eq:main_s1}
	 \|S_\delta f\|_2 \le C_\delta \|f\|_{p(\delta)}\,,
\end{equation}
	and function $\delta \mapsto C_\delta$ is bounded on each $[0, \Delta], \Delta < 1/2$.
\item[S2.](subcritical estimate around $\delta=1/2$) For any $p>1$ and sufficiently small $\theta$
(in particular, $p>1+\theta^2$) there is $C$ such that for all
	$\delta \in [(1-\theta)/2, 1/2]$, $n\ge1$ and functions $f$ we have
\begin{equation}\label{eq:main_s2}
	 \|S_\delta f\|_2 \le C \|f\|_p 
\end{equation}
\end{enumerate}

First we show S1. In accordance with~\eqref{eq:kr_sp}
\begin{equation}\label{eq:p_2}
	 \|S_\delta f\|_2^2 = \sum_{a=0}^n \left|K_{\delta n}(a)\over K_{\delta n}(0)\right|^2 \|f_a\|_2^2\,, 
\end{equation}
where we denoted
	$$ f_a \eqdef \Pi_a f\,.$$

\begin{figure}[t]
\centering
\includegraphics[width=.4\textwidth]{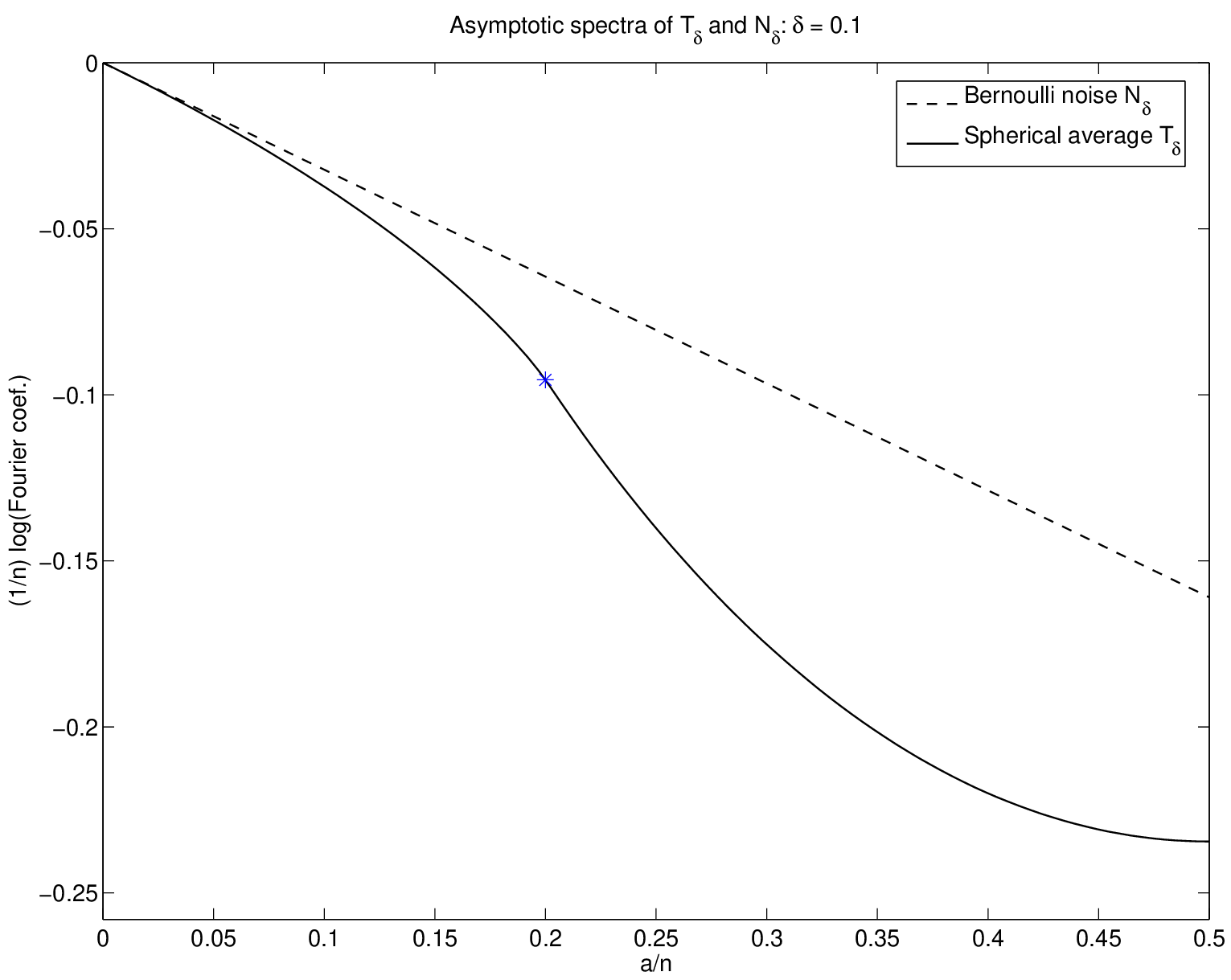}\hskip30pt
\includegraphics[width=.4\textwidth]{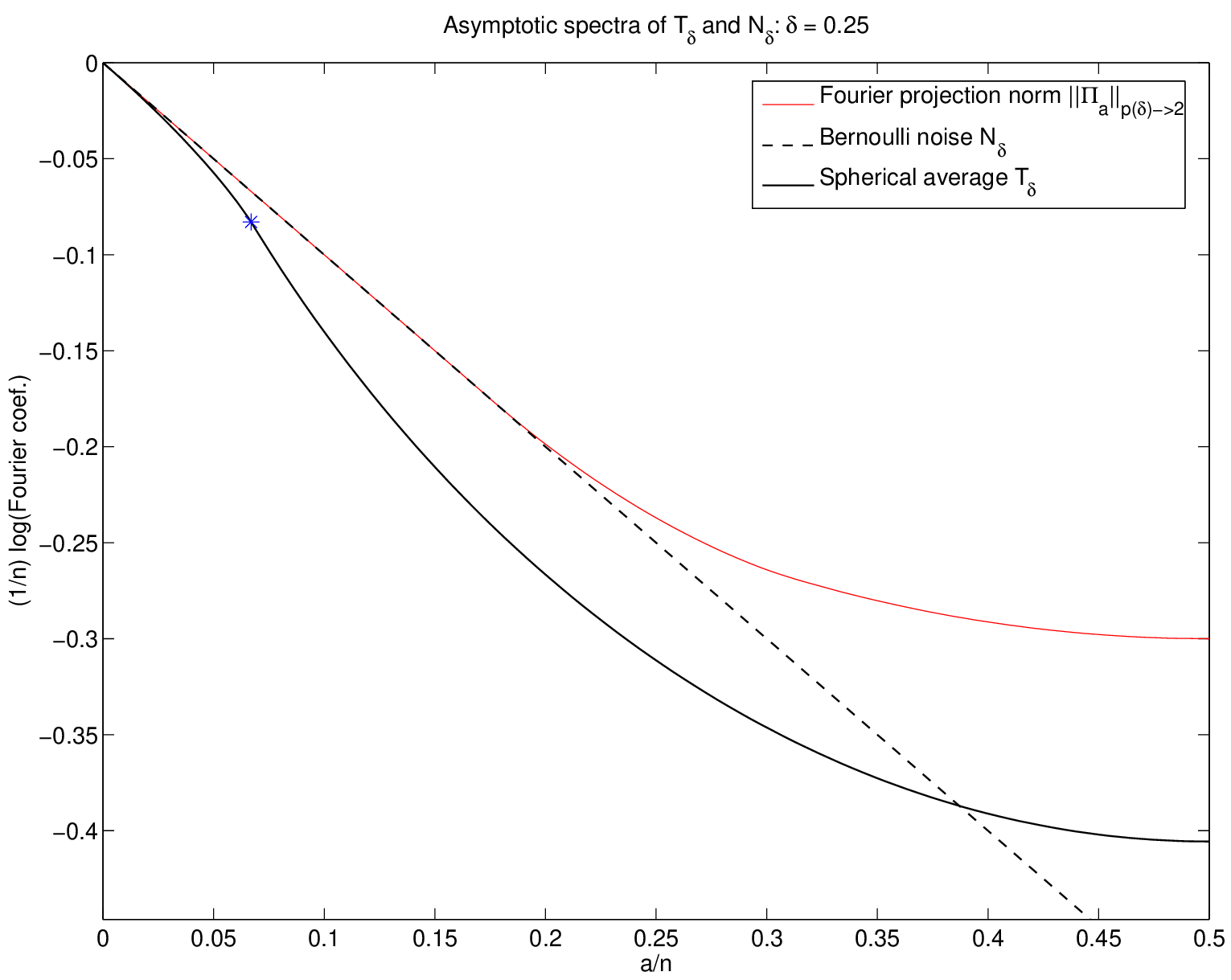}
\caption{Comparison of exponents of $a$-th eigenvalue of $S_\delta$ and $N_\delta$. For larger
$\delta$ we also show the negative of the exponent of $\|\Pi_a\|_{p(\delta)\to 2}$,
$p(\delta)=1+(1-2\delta)^2$. As before asterisks denote the critical value $\xi_{crit}(\delta)$,
i.e. the smallest root of Krawtchouk polynomial $K_{\delta n}(\cdot)$.} \label{fig:compare_expon}
\end{figure}
The scheme of our proof is illustrated by Fig.~\ref{fig:compare_expon}:
\begin{enumerate}
\item First, we show that summation in~\eqref{eq:p_2} can be truncated to $a\le {n\over 2}$.
\item Second, we show that for small values of $\delta$ eigenvalues of $S_\delta$ are upper-bounded
by a constant multiple of eigenvalues of $N_\delta$ defined in~\eqref{eq:ndelta}. This is the
content of Lemma~\ref{th:rzz}.
\item Third, for larger values of $\delta$ we show that although eigenvalues of $S_\delta$ can be
exponentially larger than those of $N_\delta$, such eigenvalues correspond to large $a$ for which
${\|f_a\|_2\over \|f\|_p}$ is exponentially smaller. 
\end{enumerate}

For the first step note that any $f$ can be written as
	$$ f = f_{even} + f_{odd}\,,$$
where each of the summands is supported on vectors $x\in\FF_2^n$ of even/odd weight. Note that
$S_\delta f_{even}$ and $S_\delta f_{odd}$ are also of opposite parity. Thus,
	$$ \|S_\delta f\|_2^2 = \|S_\delta f_{even}\|_2^2 + \|S_{\delta} f_{odd}\|_2^2\,.$$
On the other hand, we have 
\begin{align} \left(\|f_{even}\|_p^2 + \|f_{odd}\|_p^2\right)^{1\over 2} &\le 
		\left\| \sqrt{f_{even}^2 + f_{odd}^2}\right\|_p\label{eq:p_2a}\\
		&=\|f\|_p\,,\label{eq:p_2b}
\end{align}
where~\eqref{eq:p_2a} is from Minkowski's inequality and~\eqref{eq:p_2b} is because the supports of
$f_{even}$ and $f_{odd}$ are disjoint. 
Thus, if~\eqref{eq:main_s1} is established for both odd and even functions then~\eqref{eq:main_s1} follows
for all functions with the same constant $C$.

Note that for both odd and even functions we have
	$$ |\hat f(\omega)| = |\pm \hat f(\bar \omega)| = |\hat f(\bar\omega)|\,.$$
and for any such $f$ from~\eqref{eq:p_2} and~\eqref{eq:kr_sym} we get
\begin{equation}\label{eq:main1}
	 \|S_\delta f\|_2^2 \le 2 \sum_{0 \le a \le n/2} \left|K_{\delta n}(a)\over K_{\delta n}(0)\right|^2 \|f_a\|_2^2\,.
\end{equation}
\apxonly{(Inequality because of possible double-counting of $a=n/2$.)}

In the remaining we show that~\eqref{eq:main1} is upper-bounded by $C \|f\|_{p(\delta)}$ uniformly
in $f$ and $\delta \le \Delta < 1/2$. For all $\delta \in [0, \delta_0]$ from Lemma~\ref{th:rzz} we have
\begin{align} \|S_\delta f\|_2^2 &\le 2  C_1^2 \sum_{0 \le a \le n/2} (1-2\delta)^{2a} \|f_a\|_2^2\\
				&= 2 C_1^2 \|N_\delta f\|_2^2 \\
				&\le 2 C_1^2 \|f\|_{p(\delta)}^2\,,\label{eq:p_4}
\end{align}
where the last step follows from Bonami-Gross~\eqref{eq:bonami}. For $\delta \in[\delta_0, \Delta]$
we have from Lemma~\ref{th:rdd}
\begin{equation}\label{eq:p_5}
	 \left|K_{\delta n}(a)\over K_{\delta n}(0)\right| \le (1-2\delta)^{a}\,, \qquad 0 \le a
\le n \xi_{crit}(\delta)\,.
\end{equation}
On the other hand, for $a \in[n \xi_{crit}(\delta), n/2]$ we have the estimate given by Lemma~\ref{th:tough}.
Putting together~\eqref{eq:p_5} and~\eqref{eq:tough} we get similar to~\eqref{eq:p_4}:
\begin{align} \|S_\delta f\|_2^2 &\le 2 C_1^2\|N_\delta f\|_2^2 + 2\|f\|_{p(\delta)}^2 \sum_{a \in [n \xi_{crit}(\delta), n/2]}
			(C_1')^2 n e^{-2C_2 n} \\
			&\le 2C_1^2 \|N_\delta f\|_2^2 + 2 (C_1')^2 \|f\|_{p(\delta)}^2 \cdot n^2 e^{-2 C_2 n} \\
			&\le 2 (C_1^2+(C_1')^2 n^2 e^{-2C_2 n}) \|f\|_{p(\delta)}^2\,,\label{eq:p_6}
\end{align}
where in the last step we applied~\eqref{eq:bonami}. Since constants $C_1'$ and $C_2$ only depend on
$\delta_0$ and $\Delta$ we finish the proof of~\eqref{eq:main_s1} and of statement S1.

We proceed to statement S2. Showing~\eqref{eq:main_s2} is significantly simpler since $p>p(\delta)$ this
time. Take $\theta_1 = \sqrt{p-1} > \theta$ and $\delta_1 = {1-\theta_1\over 2}$. Let 
$$ \xi_1 \eqdef {\theta_1\over 1+\theta_1^2}(\theta_1 - \theta) $$
and assume that $\theta$ is so small that $\xi_{crit}(\delta) <\xi_1$ for all $\delta \in
[{1-\theta\over2},{1\over2}]$. Then, on one hand, for all 
$ 0 \le a \le n\xi_1 $
and all $\delta \in [ {1-\theta\over 2}, {1\over 2}]$ we have from Lemma~\ref{th:rqq}:
	$$ \left|K_j(a) \over K_j(0)\right| \le (1-2\delta_1)^a\,.$$
Thus, from~\eqref{eq:bonami} we get
\begin{equation}\label{eq:p_8}
	 \sum_{a\in[0, n\xi_1]} \left|K_j(a) \over K_j(0)\right|^2 \|f_a\|_2^2 \le
\|N_{\delta_1} f\|_2^2 \le \|f\|_p^2\,.
\end{equation}
On the other hand, for $a > n\xi_1$ we have for some $C_1, E>0$:
	\begin{equation}\label{eq:p_7}
	 \left|K_j(a) \over K_j(0)\right| \cdot {\|f_a\|_2\over \|f\|_p} \le C_1 \sqrt{n}
e^{-nE}\,, \qquad \forall a \in[n\xi_1, {n\over 2}]\,
\end{equation}
Indeed, from Lemma~\ref{th:kasymp} and~\eqref{eq:fp_l1} the exponent of the left-hand side
of~\eqref{eq:p_7} is upper-bounded by 
	$$ {1\over 2}\left(\ln 2 - h(\delta)\right) + \left({1\over p} - 1\right)
h(\xi)\,, \qquad \xi \eqdef {a\over n}, \, \delta \eqdef {j\over n}$$
since $\xi\in(\xi_{crit}(\delta),1/2]$.
The largest value is attained when $\delta = {1-\theta\over 2}$ and $\xi = \xi_1$, yielding
	$$ {1\over 2}\left(\ln 2 - h(\delta)\right) + \left({1\over p} - 1\right)
h(\xi)\le {1\over 2}\left(\ln 2 - h({1-\theta\over 2})\right) + 
	\left({1\over p} - 1\right) h\left({\theta_1(\theta_1-\theta)\over 1+\theta_1^2}\right)\,.$$
Since $p>1$ as $\theta \to 0$ the function on the right-hand side becomes negative. Thus the
exponent of left-hand side in~\eqref{eq:p_7} is negative for sufficiently small $\theta$.

Estimating the sum in~\eqref{eq:main1} via~\eqref{eq:p_8} and~\eqref{eq:p_7} we get similar
to~\eqref{eq:p_6} that 
$$ \|S_\delta f\|_2^2  \le 2 (1+(C_1)^2 n^2 e^{-2 E n}) \|f\|_p^2 \qquad \forall \delta
\in[{1-\theta\over 2}, {1\over 2}]\,.$$
This completes the proof of~\eqref{eq:main_s2} and statement S2.

We proceed to lower bounds on $\|S_\delta\|_{p\to2}$. 
To show~\eqref{eq:main_cex} consider function 
	$$ f(x) = \prod_{j=1}^n (1+\epsilon \chi_j) = \sum_{t=0}^n (1+\epsilon)^{n-t}
(1-\epsilon)^t 1_{\Sph_t} = \sum_{k=0}^n \epsilon^{k} K_k(x)\,.$$
On one hand,
\begin{align} \|f\|_p &= \left({(1+\epsilon)^p\over 2} + {(1-\epsilon)^p\over 2} \right)^{n\over p} \\
		   &= e^{n {p-1\over 2} \epsilon^2 + o(\epsilon^2)}\,,\quad \epsilon\to 0
\label{eq:p_10}
\end{align}
On the other hand, from Lemma~\ref{th:kasymp} and~\eqref{eq:kl_exact} we have
\begin{equation}\label{eq:p_9}
	 \|S_\delta f\|_2^2 = \sum_{a=0}^n e^{2n\left(E_{\delta}({a\over n}) - h(\delta) + {a\over n} \ln
\epsilon + {1\over 2} h({a\over n})\right) + o(n)}\,,
\end{equation}
where we also used
	$$ \|f_a\|_2 = \epsilon^a {n\choose a}^{1\over 2} = e^{a \ln \epsilon + n h({a\over n}) +
o(n)}\,.$$
For convenience, set $\xi = {a\over n}$. Then it is not hard to show from~\eqref{eq:ed_2} that
$$ E_\delta(\xi) - h(\delta) = \xi \ln(1-2\delta) + o(\xi)\,.$$
Then setting $\xi = \epsilon^2 (1-2\delta)^2$ we find that
	$$ E_{\delta}(\xi) - h(\delta) + \xi \ln \epsilon + {1\over 2} h(\xi) = {(1-2\delta)^2\over
2} \epsilon^2 + o(\epsilon^2)\,, \qquad \epsilon\to 0$$
Thus from~\eqref{eq:p_9} and~\eqref{eq:p_10} we get 
	$$ \liminf_{n\to\infty} {1\over n} \ln {\|S_\delta f\|_2 \over \|f\|_p} \ge {(1-2\delta)^2
- (p-1)\over 2} \epsilon^2 + o(\epsilon^2)\,.$$
Evidently, for $p< 1 + (1-2\delta)^2$ the norm $\|S_\delta\|_{p\to 2}$ grows exponentially in dimension.

Finally, estimate~\eqref{eq:main_young} follows from Young's inequality~\eqref{eq:young}:
\begin{align}
	 \|S_{1/2} f\|_2 &\le 2^n \|f\|_1  {\|1_{\Sph_{n/2}}\|_2\over |\Sph_{n/2}|}\\
				&= 2^n \cdot \left(2^{-n/2} 
		{n \choose \lfloor n/2 \rfloor}^{-1/2}\right) \|f\|_1\\
				&= (1+o(1)) \left(\pi n\over 2\right)^{1\over4} \|f\|_1
\end{align}
This upper-bound is tight as $f(x) = 1\{x=0\}$ shows.

\section*{Acknowledgement}

We are grateful to Prof. Y. Peres for a stimulating discussion.

\appendix
\section{Proof of Lemma~\ref{th:krawt_lp}}
\begin{proof} Let $j=\lfloor\delta n\rfloor$ and note that from Plancherel we have 
\begin{equation}\label{eq:klp_3}
		\|K_j\|_2 = \sqrt{|\Sph_j|} = \exp\left\{ {n\over 2} h(\delta) + O(\log n)\right\}\,.
\end{equation}	
	Consequently, we only consider $p\neq 2$ from now on.
	
	The lemma is shown by analyzing with exponential precision the expression
\begin{equation}\label{eq:klp_4}
		\| K_j\|_p^p = \sum_{a=0}^n 2^{-n} {n\choose a} |K_j(a)|^p \,,
\end{equation}
so that
	$$ nE(p,\delta) \le \ln \| K_j\|_p^p \le \ln(n+1) + nE(p,\delta)\,,$$
where 
\begin{equation}\label{eq:klp_4a}
		E(p,\delta) \eqdef {1\over n} \max_{a\le n/2} \ln {n\choose a} - n \ln 2 + p \ln|K_j(a)|\,,
\end{equation}	
and we used the symmetry to restrict analysis to $a\le n/2$. We will show below that for $p> 2$ the term exponentially dominating this
sum occurs at $a \le n \xi_{crit}(j/n)$, while for $p<2$ the dominating term is at $a=n/2$. 

First, consider $p>2$. From Lemma~\ref{th:kasymp}, we have 
\begin{equation}\label{eq:klp_5}
	E(p,\delta) \le \max_{0\le \xi\le 1/2} h(\xi) - \ln 2 + p E_\delta(\xi) + O\left(\log n \over n\right)\,.
\end{equation}
In the regime $\xi_{crit}(\delta) \le \xi \le 1/2$ we have 
$$ h(\xi) - \ln 2 + p E_\delta(\xi) = {1\over 2} (h(\delta)-\ln 2) + (1-p/2) h(\xi)\,,$$
which is decreasing in $\xi$, and hence we may restrict maximization in~\eqref{eq:klp_5} to $\xi \le
\xi_{crit}(\delta)$. We introduce parametrization $\xi = \xi(\omega)$ as in~\eqref{eq:ed_2}, with
$$ {\delta\over 1-\delta} \le \omega \le \sqrt{\delta\over 1-\delta}\,. $$
Then using identity
\begin{equation}\label{eq:klp_6}
	{d\over d\omega} \phi(\xi(\omega),\omega) = \xi'(\omega) \ln {1-\omega\over 1+\omega}
\end{equation}
we get that derivative of the expression under the $\max$ in~\eqref{eq:klp_5} is
\begin{equation}\label{eq:klp_7}
	{d\over d\omega}(\cdots) = \xi'(\omega) \left(\ln {1-\xi\over \xi} + p \ln {1-\omega\over 1+\omega}\right)\,.
\end{equation}
It is clear that this function is strictly increasing as $\omega$ ranges in~\eqref{eq:klp_6}.
For the right endpoint in~\eqref{eq:klp_6} we have $\xi=0$ and thus the derivative tends to $-\infty$, for the
left endpoint, notice that when $p=2$ and $\omega=\sqrt{\delta\over 1-\delta}$ the expression~\eqref{eq:klp_7} is
exactly zero and thus $>0$ for $p>2$. So there does exist a unique $\omega^*(p,\delta)$ such that~\eqref{eq:klp_7}
equals zero. Instead of finding the function $\omega^*(p,\delta)$ and $\xi^*=\xi(\omega^*)$ 
we fix an arbitrary value $\omega \in [0,1]$ and
find the $\delta$ for which $\omega^*(p,\delta)=\omega$. This gives expression for $\delta = \delta(\omega)$ given 
in~\eqref{eq:klp_b}. Plugging the values $\delta =
\delta(\omega)$ and $\xi^* = \xi^*(\delta(\omega),\omega)$ into~\eqref{eq:klp_5} we conclude that
$$ E(p,\delta) \le h(\xi^*) - \ln 2 + E_\delta(\xi^*) + O\left(\log n \over n\right),$$
where furthermore $E_\delta(\xi^*) = \phi(\xi, \omega)$. This completes proof of the upper bound in~\eqref{eq:klp}. 

To prove a matching lower bound, notice that for any fixed $\delta$ we have argued that $\omega = \sqrt{\delta\over
1-\delta}$ yields a positive value of~\eqref{eq:klp_7}. Consequently, the optimal value of $\xi^*$ in~\eqref{eq:klp_5} is
always $<\xi_{crit}(\delta)-\epsilon$ for some $\epsilon = \epsilon(p,\delta)>0$. Thus, taking $a=\lfloor \xi^*
n\rfloor$, we can apply the result of~\cite[Section IV]{KL95} establishing
$$ K_j(a) = \exp\{n E_\delta(\xi^*) + O(\log n)\}\,,$$
which shows that $E(p,\delta) \ge h(\xi^*) - \ln 2 + p E_\delta(\xi^*)   + O\left(\log n \over n\right) $ matching the
previous upper bound.

We now prove~\eqref{eq:klp_2}. The upper bound follows from $\|K_j\|_p \le \|K_j\|_2$ and~\eqref{eq:klp_3}. For the
lower bound, assume $j$ and $n$ are even. From~\eqref{eq:kr_sym} we have $K_k(n/2)=0$ for any odd $k$, and thus
from~\eqref{eq:kr_exch}, we have that roots of $K_{n/2}(\cdot)$ are precisely all odd integers in $[n]$, so that
$$ K_{n/2}(x) = c \prod_{m=1}^{n/2} (x-2m-1)\,,$$
where constant $c$ is found from $K_{n/2}(0)={n\choose n/2}$. Applying~\eqref{eq:kr_exch} again, we find
$$ K_j(n/2) = {{n\choose j} \over {n \choose n/2}} K_{n/2}(j)\,.$$
When $j$ is even, $K_j(n/2)$ is non-zero, so analyzing this similar to proof of Stirling formula we get 
$$ K_j(n/2) = \exp\{n h(\delta)/2 + O(\log n)\}\,.$$
The lower bound in~\eqref{eq:klp_2} then follows from, cf.~\eqref{eq:klp_4},
$$ 
	\| K_j\|_p^p \ge 2^{-n} {n\choose n/2} |K_j(n/2)|^p \,.
$$
\end{proof}

\ifmapx
\section{Other notes}

\subsection{Gaussianity of $K_j$}

An interesting question is the following: Does $K_j(X)$ approach gaussian distribution for $X$
uniform on $\FF_2^n$ and $n\to\infty$?

The answer turns out to be negative for all $j\ge 1$. Indeed, on one hand we can show 
	$$ \sup_n {\|K_j\|_6^6 \over \|K_j\|_2^3} < \infty $$ 
and hence the sequence $\{ {\|K_j\|_4^4\over \|K_j\|_2^4}, n=1,\ldots\}$ is uniformly integrable.
Thus,
	$$ {K_j(X) \over \|K_j\|_2} \stackrel{d}{\to} \matn(0,1) \quad \implies \quad 
			{\|K_j\|_4^4 \over \|K_j\|_2^4} \to 3\,.$$
But on the other hand, we can compute:
	$$ {\|K_j\|_4^4 \over \|K_j\|_2^4} \to \sum_{t=0}^j {j\choose t}^2 {2t \choose t}\,.$$
And the sequence of numbers on the right is:
	$$ \begin{pmatrix} j=& 0 & 1 & 2 & 3 & 4 & 5 \\
			& 1 & 3 & 15 & 93 & 639 &4653 \end{pmatrix} $$
Thus, for $j\neq 1$ the distribution of $K_j(X)$ is significantly non-gaussian.

\subsection{Open questions}

TODO:
\begin{enumerate}
\item Semenov-Shneiberg: relate to averages of symmetric functions over codes with bounded
mindist/maxdist.
\item $\|\Pi_a\|_{p\to q}$ via contour integration, Minkowski and Bonami-Gross.
\item More generally, identities and properties of $BSC(z)$ complex-noise channels.
\item Implications in terms of additive energy $E(A,B)$ when $B$ is permutation invariant?
\item Perhaps prove the following general property of groups:
	$$ \sup_x \left| {K_j(x)\over K_j(0) }\right| \cdot \|() * 2^{-n} K_j\|_{p(j)\to 2} = O(1)\,,$$
where $p(j)=1+\sigma_2^2$ determined by the spectral gap?
\item Study whether maxima of $\|f_a\|_p\over \|f_a\|_2$ and $\|f_a\|_p\over \|f\|_2$ is achieved by
permutation invariant functions when $p=2m$. (Not true for $a=1$ and $2<p<3$, but maybe for
$p=2m$?).
\end{enumerate}

\subsection{Semenov-Shneiberg and hypercontractivity of small perturbations of $\EE$}

So, using the main theorem and Semenov-Shneiberg, we can show that under mild conditions all permutation invariant noise
distributions $P_Z$ on $\FF_2^n$ such that $\|P_Z - U_Z\|_{TV}\le \epsilon$ are hypercontractive where $\epsilon$ does
not depend on dimension. Namely:
	$$ \|f*P_Z\|_p \le \|f\|_{q} $$
	for some $p<q$ and all dimensions $n$. 

However, such result is in fact trivial. Indeed, consider any noise distribution (on any finite group in fact) $P_Z$
such that 
$$ \|P_Z - U_Z\|_{TV} < {1\over 2} $$
where $U_Z$ is the Haar distribution (uniform). Then, for any $\delta_v$ we have
$$ \|P_Z - \delta_v * P_Z\|_{TV} \le 1-\delta_0 < 1 \,,$$
and thus by Dobrushin's characterization we have $\eta_{TV} < 1$ and by general results:
$$ D(P*P_Z \| Q*P_Z) \le \eta_{TV} D(P||Q) \qquad \forall P,Q$$
and by Ahlswede-Gacs we also have
$$ \|f(X)\|_p \le \|f(X+Z)\|_{q} \qquad \forall 1\le \eta_{TV} \cdot p \le q \le p < \infty \quad \forall P_X $$
This is much stronger result since it shows explicitly that a lot of noise distributions are simultaneously
hyper-contractive for ALL $P_X$.
\fi


\begin{thebibliography}{10}

\bibitem{AC76}
{\sc R.~Ahlswede and P.~Gacs}, {\em Spreading of sets in product spaces and
  hypercontraction of the {M}arkov operator}, Ann. Probab.,  (1976),
  pp.~925--939.

\bibitem{benabbas2012isoperimetric}
{\sc S.~Benabbas, H.~Hatami, and A.~Magen}, {\em An isoperimetric inequality
  for the hamming cube with applications for integrality gaps in degree-bounded
  graphs}, Unpublished, 1 (2012), p.~1.

\bibitem{AB70}
{\sc A.~Bonami}, {\em {\'E}tude des coefficients de {F}ourier des fonctions de
  $l_p(g)$}, Ann. Inst. Fourier (Grenoble), 20 (1970), pp.~335--402.

\bibitem{CB82}
{\sc C.~Borell}, {\em Positivity improving operators and hypercontractivity},
  Math. Zeit., 180 (1982), pp.~225--234.

\bibitem{PD73}
{\sc P.~Delsarte}, {\em An algebraic approach to the association schemes of
  coding theory}, Philips Research Rep. Supp.,  (1973), p.~103.

\bibitem{DSC96}
{\sc P.~Diaconis and L.~Saloff-Coste}, {\em Logarithmic {S}obolev inequalities
  for finite {M}arkov chains}, Ann. Appl. Probab., 6 (1996), pp.~695--750.

\bibitem{DS58}
{\sc N.~Dunford and J.~Schwartz}, {\em Linear Operators: General theory},
  vol.~1, Interscience Publishers, New York, 1958.

\bibitem{FS72}
{\sc C.~Fefferman and H.~S. Shapiro}, {\em A planar face on the unit sphere of
  the multiplier space $m_p, 1<p<\infty$}, Proc. AMS, 36 (1972).

\bibitem{frankl1987forbidden}
{\sc P.~Frankl and V.~R{\"o}dl}, {\em Forbidden intersections}, Trans. Amer.
  Math. Soc., 300 (1987), pp.~259--286.

\bibitem{RG68}
{\sc R.~G. Gallager}, {\em Information Theory and Reliable Communication},
  Wiley, New York, 1968.

\bibitem{GMPT10}
{\sc K.~Georgiou, A.~Magen, T.~Pitassi, and I.~Tourlakis}, {\em Integrality
  gaps of 2-o(1) for vertex cover sdps in the lov{\'a}sz--schrijver hierarchy},
  SIAM Journal on Computing, 39 (2010), pp.~3553--3570.

\bibitem{LG75}
{\sc L.~Gross}, {\em Logarithmic sobolev inequalities}, Amer. J. Math., 97
  (1975), pp.~1061--1083.

\bibitem{IS98}
{\sc M.~E.~H. Ismail and P.~Simeonov}, {\em Strong asymptotics for {K}rawtchouk
  polynomials}, J. Comp. and Appl. Math., 100 (1998), pp.~121--144.

\bibitem{KKL88}
{\sc J.~Kahn, G.~Kalai, and N.~Linial}, {\em The influence of variables on
  {B}oolean functions}, in Proc. 29th Ann. Symp. on Foundations of Comp. Sci.,
  Los Alamitos, CA, 1988, pp.~68--80.

\bibitem{KL95}
{\sc G.~Kalai and N.~Linial}, {\em On the distance distribution of codes},
  {IEEE} Trans. Inf. Theory, 41 (1995), pp.~1467--1472.

\bibitem{kauers2014hypercontractive}
{\sc M.~Kauers, R.~O'Donnell, L.-Y. Tan, and Y.~Zhou}, {\em Hypercontractive
  inequalities via sos, and the frankl--r{\"o}dl graph}, in Proceedings of the
  Twenty-Fifth Annual ACM-SIAM Symposium on Discrete Algorithms, SIAM, 2014,
  pp.~1644--1658.

\bibitem{kneser1952theoreme}
{\sc H.~Kneser}, {\em Sur un th{\'e}oreme fondamental de la th{\'e}orie des
  jeux}, Comptes Rendus Acad. Sci. Paris, 234 (1952), pp.~2418--2420.

\bibitem{KL01}
{\sc I.~Krasikov and S.~Litsyn}, {\em Survey of binary {K}rawtchouk
  polynomials}, DIMACS series: Codes and association schemes, 56 (2001),
  pp.~199--212.

\bibitem{MRRW77}
{\sc R.~McEliece, E.~Rodemich, H.~Rumsey, and L.~Welch}, {\em New upper bounds
  on the rate of a code via the {D}elsarte-{M}ac{W}illiams inequalities},
  {IEEE} Trans. Inf. Theory, 23 (1977), pp.~157--166.

\bibitem{LM97}
{\sc L.~Miclo}, {\em Remarques sur l'hypercontractivit{\'e} et l'{\'e}volution
  de l'entropie pour des cha{\^\i}nes de {M}arkov finies}, in S{\'e}minaire de
  {P}robabilit{\'e}s {XXXI}, Springer, 1997, pp.~136--167.

\bibitem{MORSS-nicd}
{\sc E.~Mossel, R.~O'Donnell, O.~Regev, J.~E. Steif, and B.~Sudakov}, {\em
  Non-interactive correlation distillation, inhomogeneous markov chains, and
  the reverse bonami-beckner inequality}, Israel Journal of Mathematics, 154
  (2006), pp.~299--336.

\bibitem{odonnell-book}
{\sc R.~O'Donnell}, {\em Analysis of boolean functions}, Cambridge University
  Press, 2014.

\bibitem{YP12-isit}
{\sc Y.~Polyanskiy}, {\em Hypothesis testing via a comparator}, in Proc. 2012
  IEEE Int. Symp. Inf. Theory (ISIT), Cambridge, MA, July 2012.

\bibitem{YP13-htstruct}
\leavevmode\vrule height 2pt depth -1.6pt width 23pt, {\em Hypothesis testing
  via a comparator and hypercontractivity}, preprint,  (2013).

\bibitem{IS70}
{\sc I.~Segal}, {\em Construction of non-linear local quantum processes: {I}},
  Ann. Math., 92 (1970), pp.~462--481.

\bibitem{SS88}
{\sc E.~M. Semenov and I.~Y. Shneiberg}, {\em Hypercontractive operators and
  {K}hinchin's inequality}, Func. Analysis and Appl., 22 (1988), pp.~244--246.

\end{thebibliography}
\end{document}